\documentclass[12pt]{article}
\usepackage{babel}[english]
\usepackage{amsmath,amsfonts,amsthm,amssymb}
\usepackage[dvips]{graphics}
\usepackage{epsfig,indentfirst,relsize,comment}
\usepackage[left=2cm,right=2cm,top=3.5cm,bottom=3.5cm]{geometry}
\usepackage[dvipsnames]{xcolor}
\allowdisplaybreaks

\newtheoremstyle{theorem}
  {10pt}
  {10pt}
  {\sl}
 {}
  {\bf}
  {.\,}
  { }
  {}
\theoremstyle{theorem}
\newtheorem{theorem}{Theorem}[section]
\newtheorem{corollary}{Corollary}[section]
\newtheorem{definition}{Definition}[section]

\newtheorem{remark}{Remark}[section]

\numberwithin{equation}{section}

\newtheoremstyle{defi}
{10pt}
{10pt}
{\rm}
{}
{\bf}
{.\,}
{ }
{}
\theoremstyle{defi}

\newcommand{\D}{\partial}
\renewcommand{\div}{\textrm{div}\,}
\newcommand{\R}{\mathbb{R}}
\newcommand{\bee}{\mathbf{u}^E}
\newcommand{\bv}{\mathbf{v}}

\definecolor{mypink}{RGB}{219, 48, 122}
\definecolor{myblue}{RGB}{0, 0, 122}

\begin{document}
\baselineskip = 13.5pt

\title{Vanishing viscosity limit for the compressible Navier-Stokes 
equations with non-linear density dependent viscosities}

\author{Luca Bisconti$^{1}$ \,\, \ Matteo Caggio$^{2\,*}$ \,\, \ Filippo Dell'Oro$^{3}$
  \\
  \\
  {\small  1. Università degli Studi di Firenze - Dipartimento di Matematica e Informatica “U. Dini”}\\\smallskip
  {\small Viale Morgagni 67/a, I-50134 Firenze, Italy}\\
  {\small  2. Institute of Mathematics of the Academy of Sciences of the Czech Republic} \\\smallskip
  {\small \v Zitn\' a 25, 11567 Praha 1, Czech Republic}\\
  {\small  3. Politecnico di Milano - Dipartimento di Matematica}\\
  {\small Via Bonardi 9, 20133 Milano, Italy}\\
  {\small *\,corresponding author:\ caggio@math.cas.cz}
\date{}
}

\maketitle

\begin{abstract}
In a three-dimensional bounded domain $\Omega$ we consider the
compressible Navier-Stokes equations for a barotropic fluid with general
non-linear density dependent viscosities and no-slip boundary conditions. 
A nonlinear drag term is added to the momentum equation.
We establish two conditional Kato-type criteria for the convergence of
the weak solutions to such a system towards the strong solution of
the compressible Euler system when the viscosity coefficient
and the drag term parameter tend to zero.
\end{abstract}

\noindent
{\small{\bf Keywords:} compressible Navier-Stokes equations, density dependent viscosities,
vanishing viscosity limit, boundary layer.}\smallskip

\noindent
{\small{\bf 2010 Mathematics Subject Classifications:} 35Q30, 35Q35, 76N10}

\section{Introduction}\setcounter{equation}{0}
\noindent
Let $\Omega \subset \mathbb{R}^3$ be a bounded domain with
sufficiently smooth boundary $\partial\Omega$, and let $T>0$ be
arbitrarily fixed. We consider the compressible Navier-Stokes
equations describing the motion of a barotropic fluid with density
dependent viscosities that may vanish on vacuum
\begin{gather} \label{cont}
  \partial_{t}\varrho +\textrm{div}(\varrho \mathbf{u})=0,\\\noalign{\vskip0.5mm}
  \label{mom}
  \partial_{t} (\varrho \mathbf{u})+\textrm{div}(\varrho\mathbf{u}\otimes\mathbf{u}) +
  \nabla p(\varrho) -\varepsilon [\hspace{0.3mm}\textrm{div}
  (2\mu(\varrho) \mathbb{D}(\mathbf{u})) +\nabla (\lambda(\varrho)
  \text{div} \mathbf{u})] + r_1 \varrho \hspace{0.3mm}
  |\mathbf{u}|\mathbf{u} = 0,
\end{gather}
where $\varepsilon$ and $r_1$ are positive parameters. The unknown
variables
\begin{equation*}
 \varrho : (0,T)\times \Omega \to \R\qquad \text{and}
\qquad \mathbf{u} :(0,T)\times \Omega \to \R^3,
\end{equation*}
with $\varrho=\varrho\left(t,x\right)$, $\mathbf{u}=\mathbf{u}\left(t,x\right)$,
represent the density and the velocity vector field,
respectively.  The term $p=p(\varrho)$ represent the pressure of the
fluid and is given by a standard power law
\begin{equation*}
 p(\varrho) = \varrho^\gamma, \quad\,\, \gamma > 1.
\end{equation*}
Moreover
$\mathbb{D}(\mathbf{u})=\frac12(\nabla \mathbf{u} + \nabla^{\top}
\mathbf{u})$ is the symmetric part of $\nabla \mathbf{u}$.  The
density dependent viscosity coefficients $\mu(\varrho)$ and
$\lambda(\varrho)$ satisfy the following algebraic relation
\begin{equation} \label{lm} \lambda(\varrho) = 2(\varrho\mu'(\varrho)
  -\mu(\varrho))
\end{equation}
and further structural assumptions specified later. The last term in the
momentum equation \eqref{mom} corresponds to a turbulent drag force
(see, e.g. \cite{BDGV}).

The system is supplemented with the no-slip boundary condition
\begin{equation} \label{bc}
  \varrho \mathbf{u}|_{\partial \Omega}=0.
\end{equation}
Moreover, as in \cite{BDGV}, we prescribe the additional boundary
condition
\begin{equation} \label{bc-bis}
  \mu(\varrho)\nabla s(\varrho) \times \mathbf{n} |_{\partial  \Omega}= 0,
\end{equation}
where $\mathbf{n}$ is the unit vector normal to $\partial\Omega$, and
$s'(\varrho)=\mu'(\varrho)/\varrho$.
The latter expresses the fact
that the density should be constant on the connected components of
$\partial\Omega$ (see \cite{BDGV}). Finally, the initial state of the system is prescribed
\begin{equation*}
  \varrho (0,\cdot) =
  \varrho_{0}, \qquad \varrho \mathbf{u} (0,\cdot) = \varrho_{0}\mathbf{u}_{0}.
\end{equation*}
Formally, letting $\varepsilon,r_1\rightarrow0$, one obtains the
compressible Euler system
\begin{gather} \label{cont-E}
  \partial_t\varrho^E +\textrm{div} (\varrho^E \mathbf{u}^E) = 0,	\\
  \label{mom-E}
  \partial_t (\varrho^E \mathbf{u}^E) + \div (\varrho^E \mathbf{u}^E
  \otimes \mathbf{u}^E) + \nabla p(\varrho^E) =0,
\end{gather}
with the boundary condition
\begin{equation} \label{bc-E}
  \bee\cdot\mathbf{n}|_{\partial \Omega}=0
\end{equation}
and the initial conditions
\begin{equation*}
\varrho^E(0)=\varrho_0^E, \qquad
 \varrho^E \mathbf{u}^E (0) = \varrho_0^E \mathbf{u}_0^E.
\end{equation*}
As is well-known, a rigorous justification of such a singular limit on
domains with non-empty boundary is a challenging problem, due to the
appearance of boundary layers.  A remarkable approach was introduced
in the seminal paper of Kato~\cite{Ka}, where the author studied the
vanishing viscosity limit (also called inviscid limit) of the
incompressible Navier-Stokes equations with no-slip boundary
conditions. He proved that, if the energy dissipation rate of the
viscous flow in a boundary layer of width proportional to the
viscosity vanishes, then the solutions of the incompressible
Navier–Stokes equations converge to the solutions of the
incompressible Euler equations. The analysis requires the construction
of an artificial ``fake" boundary layer, needed to compensate the
discrepancy between the boundary conditions of the original system and
the boundary conditions of the target one.  The result of~\cite{Ka}
has been improved and generalized by several authors. Without any
claim of completeness, we mention \cite{Ke,Wa} among many others.  We
also mention the work of Constantin et al.~\cite{Co} where the
inviscid limit has been performed under the so-called Oleinik
condition of no back-flow in the trace of the Euler flow and a lower
bound for the Navier-Stokes vorticity in a Kato-like boundary layer.

For the compressible Navier-Stokes equations with no-slip boundary
conditions and constant viscosity coefficients, that is, when the
viscous term
\begin{equation*}
  -\varepsilon [\hspace{0.3mm}\textrm{div} (2\mu(\varrho)
\mathbb{D}(\mathbf{u})) +\nabla (\lambda(\varrho) \text{div}
\mathbf{u})]
\end{equation*}
is replaced by
$-\varepsilon\, \textrm{div}\, \mathbb{S}(\nabla \mathbf{u})$, where
\begin{equation*}
  \mathbb{S}(\nabla \mathbf{u}) = \mu\Big( \nabla \mathbf{u} + 
  \nabla ^\top \mathbf{u} -\frac{2}{3}\big( \text{div} \mathbf{u}\,\mathbb{I}
  \big) \Big) +\eta(\text{div}\mathbf{u})\mathbb{I},
\end{equation*}
for some constants $\mu>0$ and $\eta\geq0$, some results
have been obtained (without the presence of the drag term in the momentum equation)
by Sueur \cite{Su}.  He established the following sufficient
condition for the vanishing viscosity limit to hold:
\begin{equation} \label{Sueur-cond}
  \int_{0}^{T} \int_{\Gamma_{c\varepsilon}} \left(
    \varepsilon\,\frac{\varrho|\mathbf{u}|^2}{d_\Omega^2}
    +\varepsilon\,\frac{\varrho^2(\mathbf{u}\cdot
      \mathbf{n})^2}{d_\Omega^2} +\varepsilon\,|\mathbb{S}(\nabla
    \mathbf{u})|^2 \right)dxdt \to 0 \,\,\mbox{ as }\,\, \varepsilon
  \to 0.
\end{equation}
Here and in the sequel, we denote by $d_\Omega = d_\Omega (x)=d(x, \D\Omega)$ the distance of
$x \in \Omega$ to the boundary $\partial \Omega$, and by
\begin{equation*}
\Gamma_{c\varepsilon} = \{ x\in \Omega :  d_\Omega(x) < c  \varepsilon \},
\end{equation*}
for a constant $c > 0$ (in the case of \eqref{Sueur-cond}, the
constant $c>0$ is arbitrarily fixed).  Condition~\eqref{Sueur-cond}
can be regarded as a generalization of the result of Kato in the
context of compressible fluids. Indeed, in the light of the Hardy
inequality, when the density is constant \eqref{Sueur-cond} is implied by
\begin{equation*}
\varepsilon \int_{0}^{T} \int_{\Gamma_{c\varepsilon}}
    |\nabla \mathbf{u}|^2
    dxdt \to 0 \,\,\mbox{ as }\,\, \varepsilon \to 0,
  \end{equation*}
  which is exactly the condition employed by Kato in~\cite{Ka}.
Subsequently, Bardos and Nguyen \cite{BaNg} proved
the convergence to the so-called dissipative solutions of the compressible Euler equations
under the alternative condition
\begin{equation} \label{Bardos-Nguyen-cond}
   \int_{0}^{T} \int_{\Gamma_{c\varepsilon}}
    \left(
    \varrho^\gamma
    +\varepsilon\,\frac{\varrho|\mathbf{u}|^2}{d_\Omega^2}
    +\varepsilon\, |\nabla \mathbf{u}|^2
    \right)dxdt \to 0\,\, \mbox{ as }\,\, \varepsilon \to 0.
  \end{equation} 
  The achievements of \cite{Su} have been somehow improved by Wang and
  Zhu \cite{WaZh}, where the authors assumed sufficient conditions
  only on the tangential or the normal component of velocity, at the
  cost of increasing the width of the boundary layer. See also the
  recent article \cite{cinesiJMFM} which can be regarded as the
  analogous of \cite{Co} for compressible fluids.

  We now come to the literature concerning the vanishing viscosity
  limit for the compressible Navier-Stokes system with density
  dependent viscosity coefficients. In the three-dimensional torus,
  the convergence has been analyzed by Bresch et.\ al.\ \cite{BNV-2}
  employing the so-called ``augmented version" of the compressible
  Navier-Stokes system.  Further results have been obtained by Geng et
  al.~\cite{Ge,Ge-2} in the whole space $\mathbb{R}^3$ and for viscosity
  coefficients of the form $\varrho^\beta$ with $\beta > 1$. We also
  mention the articles~\cite{ChChZh,DiZhu} which study the case when the
  viscosity coefficients obey a lower order power law of the density,
  namely, $\varrho^\beta$ with $0 < \beta \leq 1$. The case of bounded
  domains with non-empty boundary has been recently treated in
  \cite{BC} for the particular choice
\begin{equation*}
\mu(\varrho)=\varrho\qquad \text{and}\qquad \lambda(\varrho)=0,
\end{equation*}
 namely, for linear density dependent viscosity coefficients. In this article,
the passage to the limit has been performed under the assumption:
\begin{equation} \label{BiCa-cond}
\int_{0}^{T} \int_{\Gamma_{c\varepsilon}}
\left(\frac{\varrho^\gamma}{\varepsilon}
+\varepsilon^{\frac{\gamma-1}{\gamma}} 
\frac{\varrho|\mathbf{u}|^2}{d_\Omega^2}\right) dxdt
\to 0 \,\, \textrm{ as } \,\, \varepsilon \to 0.
\end{equation}
Note that, compared to \eqref{Sueur-cond} and \eqref{Bardos-Nguyen-cond},
relation \eqref{BiCa-cond} requires a decay rate for the
$L^\gamma$-norm, $\gamma>1$, of the density and a stronger
condition on the kinetic energy term. But notably, no assumptions on
the viscous terms are imposed.  The strategy adopted in \cite{BC}
relies on the so-called relative energy estimates and, in particular,
on the introduction of a relative energy functional ``measuring" the
distance between a weak solution of the compressible Navier-Stokes
system and a strong solution of the compressible Euler system. A
relative energy inequality satisfied by the weak solutions of the
Navier-Stokes equations has been derived, which however involves some
test functions that must satisfy the no-slip boundary conditions. 
For this reason, as done in the previous
contributions mentioned above, a ``correction" based on a Kato-type
``fake" boundary layer was needed.  Still, different from the standard
construction of Feireisl et al.\ \cite{FNS-2011,FNJ-2012}, the
relative energy inequality constructed in \cite{BC} is derived from
the aforementioned ``augmented version" of the compressible
Navier-Stokes system, in the same spirit of Bresch et
al. \cite{BNV-1,BNV-2} through the use of the Bresch-Desjardins
entropy inequality (see also \cite{BDZ, BGL, CD, Cia} for
recent applications).

\subsection*{Description of our results}
\noindent
The aim of the present paper is twofold, and can be summarized in the following points.
\begin{itemize}
\item We extend the results of \cite{BC} to more general
  and possibly non-linear density dependent viscosity coefficients
  $\mu(\varrho)$ and $\lambda(\varrho)$. Actually, we also weaken the
  second assumption in \eqref{BiCa-cond} by replacing the factor
  $\varepsilon^{(\gamma-1)/\gamma}$ with $\varepsilon$. See the
  forthcoming Theorem~\ref{main-1}.

\item We establish a new conditional criterion analogous to
  \eqref{Sueur-cond} for the validity of the vanishing viscosity limit
  in the context of compressible fluids with density dependent
  viscosities. See Theorem~\ref{main-2} below.
\end{itemize}

\noindent  
In the same spirit of \cite{BC}, we employ the relative energy
inequality approach and the construction of the ``fake" boundary
layer. However, unlike \cite{BC}, we completely avoid the use of the
``augmented version" of the compressible Navier-Stokes system as well
as the use of the Bresch-Desjardins entropy inequality, following
instead the construction of Feireisl at al.\ \cite{FNS-2011, FNJ-2012}. 

\subsection*{Plan of the paper}
\noindent
In the next Section~2, we introduce the assumptions satisfied by the viscosity coefficients
and the weak formulation of the compressible Navier-Stokes system. Moreover, we discuss
the strong solution of the compressible Euler system. In the subsequent Section~3, we state
our main results Theorems~\ref{main-1}-\ref{main-2}, while in Section~4 we derive the relative
energy inequality. The remaining Sections~5-8 are devoted to the vanishing viscosity limit and
the proofs of Theorems~\ref{main-1}-\ref{main-2}.

\section{Notation and preliminaries}
\noindent
We denote by $C_c^\infty([0,T]\times \Omega;\mathbb{R}^d)$, $d\geq1$, the space 
compactly supported smooth functions with values in $\mathbb{R}^d$,
and by $L^p(\Omega)$, $p\geq 1$, the standard
Lebesgue spaces. Sobolev spaces $W^{k,p}(\Omega)$, $k\geq 1$, consist of
functions in $L^p(\Omega)$ with weak derivatives of order $k$ that still
belong to $L^p(\Omega)$. The Bochner evolution spaces, for time-dependent
functions with values in Banach spaces $X$ are denoted by $L^p(0,T;X)$ and
$W^{k,p}(0,T;X)$. Also, $C_w(0,T;X)$ denotes the space of weakly continuous
functions.

\subsection{Assumptions on viscosity coefficients $\mu$ and $\lambda$}
\noindent 
In the same spirit of \cite{BDGV,MeVa}, and recalling that $\gamma>1$, in addition to the key relation \eqref{lm}
we assume that there exists a positive constant
\begin{equation}
\label{constraint-nu}
\nu \in \left[\frac{1}{3\gamma-2},1\right)\! ,
\end{equation}
such that the following conditions on the viscosity coefficients are satisfied:
\begin{align} 
\label{ass-1}
&\mu'(\varrho) \geq \nu, \quad \mu(0) \geq 0,\\\noalign{\vskip1mm}
\label{ass-2}
&|\lambda'(\varrho)| \leq \frac{1}{\nu} \mu'(\varrho),\\
\label{ass-3}
&\nu \mu(\varrho) \leq 2\mu(\varrho) + 3 \lambda(\varrho) \leq \frac{1}{\nu} \mu(\varrho).
\end{align}
When $\gamma \geq 3$, we also require that
\begin{equation}
\label{ass-4}
 \liminf_{\varrho \to \infty} \frac{\mu(\varrho)}{\varrho^{\gamma/3+\vartheta}} > 0
\end{equation}
for some small $\vartheta > 0$.
As remarked in \cite{BDGV,MeVa}, hypothesis \eqref{lm} and \eqref{ass-3} yield the bounds
\begin{equation} \label{rho-split}
\begin{cases}
\mu(\varrho) \leq c_1\varrho^{\left(\frac{2}{3}+\frac{\nu}{3}\right)}, \quad &\text{if }\varrho\leq 1,\\
\noalign{\vskip0.7mm}
\mu(\varrho) \leq c_1\varrho^{\left(\frac{2}{3}+\frac{1}{3\nu}\right)}, \quad &\text{if }\varrho\geq 1,
\end{cases}
\end{equation}
for some $c_1>0$.
Moreover, it follows from \eqref{ass-3} that
\begin{equation} \label{C-bound}
|\lambda(\varrho)|\leq c_2 \mu(\varrho),
\end{equation}
for all $\varrho>0$ and some $c_2>0$.

\begin{remark}
The functions $\mu(\varrho)=\varrho$ and $\lambda(\varrho)=0$ satisfy
all the assumptions stated above, including~\eqref{lm}, provided that
$\gamma\geq 4/3$.
\end{remark}

\subsection{Weak solutions to the compressible Navier-Stokes system}
\noindent
We introduce the definition of the weak solution to system
\eqref{cont}--\eqref{mom}.  As is well-known, the degenerate
viscosities prevent $\mathbf{u}$ to be defined on the vacuum region
(see e.g.\ \cite{MeVa}).  As a consequence, not $\mathbf{u}$ nor
$\nabla \mathbf{u}$ are defined almost everywhere. On the other hand,
the quantities $\sqrt{\varrho}$ and $\sqrt{\varrho} \mathbf{u}$ turn
out to be well-defined (see e.g.\ \cite{MeVa}). Consequently, the
density $\varrho\geq0$ (thought as $(\sqrt{\varrho})^2$) and the
momentum $\varrho\mathbf{u}$ (thought as
$\sqrt{\varrho}\sqrt{\varrho} \mathbf{u}$) are also
well-defined. Accordingly, the correct way to interpret the
viscous term is to think it as
\begin{equation*}
  \textrm{div} (2\mu(\varrho) \mathbb{D}(\mathbf{u}))
  +\nabla (\lambda(\varrho) \text{div} \mathbf{u})=
  2 \,\textrm{div} \Big(\sqrt{\mu(\varrho)}\mathbb{S}_\mu 
  +\frac{\lambda(\varrho)}{2\mu(\varrho)}\textrm{Tr}\big(\sqrt{\mu(\varrho)}\mathbb{S}_\mu\big)\mathbb{I} \Big),
\end{equation*}
where the matrix-valued function $\mathbb{S}_\mu$ is the symmetric
part of the matrix-valued function $\mathbb{T}_\mu$ defined through
the identity (see e.g.\ \cite{BVY})
\begin{equation}
  \label{defT}
  \sqrt{\mu(\varrho)}\mathbb{T}_\mu = \nabla \Big(\sqrt{\varrho}\mathbf{u}\,
  \frac{\mu(\varrho)}{\sqrt{\varrho}} \Big)
  -\sqrt{\varrho} \mathbf{u} \otimes \sqrt{\varrho}\nabla s (\varrho),\quad\,\,\,
  \text{in $\mathcal{D}'((0,T)\times \Omega)$},
\end{equation}
with $s'(\varrho) = \mu'(\varrho)/\varrho.$

\begin{remark}
 In the sequel, especially in providing sufficient conditions for the vanishing 
 viscosity limit (in the spirit of \eqref{Sueur-cond}-to-\eqref{BiCa-cond}), 
 as a matter of notation, we write $\varrho |\mathbf{u}|^2$ in place of $|\sqrt{\varrho}\mathbf{u}|^2$.
\end{remark}
 
\begin{definition} \label{def-ws} Let $T > 0$ be arbitrarily fixed. We
  say that $(\varrho, \mathbf{u})$ with $\varrho\geq0$ is a weak
  solution of \eqref{cont}-\eqref{mom} with boundary conditions
  \eqref{bc}-\eqref{bc-bis} on the time interval $[0,T]$ if the
  following hold.
  \smallskip
  
  \noindent {\bf (1)} Regularity properties
  \begin{align*}
    \varrho \in L^\infty(0,T;L^\gamma(\Omega)), \quad \sqrt{\varrho}
    \mathbf{u},\,\,  \sqrt{\varrho}\,\nabla s(\varrho),\,\,
    \mu(\varrho)/\sqrt{\varrho}
    \in L^\infty(0,T;L^2(\Omega)).
  \end{align*}
  
  \noindent {\bf (2)} The continuity equation \eqref{cont} is
  satisfied in the following sense
  \begin{equation*}
-\int_0^\tau \int_{\Omega} \varrho \partial_t \phi\, dx dt -
\int_0^\tau \int_{\Omega} \varrho \mathbf{u} \cdot \nabla \phi \,dx
dt = \int_{\Omega} \varrho(0,\cdot) \phi(0,\cdot) \,
dx-\int_{\Omega} \varrho(\tau,\cdot) \phi(\tau,\cdot) \, dx
\end{equation*}
for all $\tau\in[0,T]$ and every
$\phi \in C_c^\infty([0,T]\times\overline{\Omega};\mathbb{R})$.
\medskip

\noindent {\bf (3)} The momentum equation \eqref{mom} is satisfied in
the following sense
\begin{gather*}
  -\int_0^\tau \int_{\Omega} \varrho \mathbf{u}\cdot \partial_t
  \boldsymbol{\varphi} \, dxdt -\int_0^\tau \int_{\Omega}
  (\sqrt{\varrho} \mathbf{u} \otimes \sqrt{\varrho} \mathbf{u}) :
  \nabla \boldsymbol{\varphi} \, dxdt -\int_0^\tau \int_{\Omega}
  p(\varrho) \textrm{div}\,\boldsymbol{\varphi} \,dxdt
  \\\noalign{\vskip0.5mm} +\varepsilon\int_0^\tau \int_{\Omega} \Big[2
  \Big(\sqrt{\mu(\varrho)}\mathbb{S}_\mu
  +\frac{\lambda(\varrho)}{2\mu(\varrho)}\textrm{Tr}\big(\sqrt{\mu(\varrho)}\mathbb{S}_\mu\big)\mathbb{I}
  \Big) : \nabla \boldsymbol{\varphi} \Big] dxdt
  +r_1\int_{0}^\tau\int_{\Omega}|\sqrt{\varrho}\mathbf{u}|\sqrt{\varrho}\mathbf{u}\cdot
  \boldsymbol{\varphi}dxdt\\\noalign{\vskip0.3mm} = \int_{\Omega}
  (\varrho\mathbf{u} \cdot \boldsymbol{\varphi})(0,\cdot) \,dx
  -\int_{\Omega} (\varrho\mathbf{u} \cdot
  \boldsymbol{\varphi})(\tau,\cdot) \,dx
\end{gather*}
for all $\tau\in[0,T]$ and every
$\boldsymbol{\varphi} \in C_c^\infty([0,T]\times\Omega;\mathbb{R}^3)$,
where the matrix-valued function $\mathbb{S}_\mu$ is the symmetric
part of the matrix-valued function $\mathbb{T}_\mu$ defined in
\eqref{defT}.
\medskip

\noindent {\bf (4)} The following energy inequality is satisfied 
\begin{equation} \label{ee}
  \int_{\Omega}\Big[\frac{1}{2}\varrho |\mathbf{u}|^{2}
  +\frac{\varrho^\gamma}{\gamma - 1}\Big](\tau,\cdot) dx
  +\varepsilon\nu\int_{0}^\tau\int_{\Omega}\left|\mathbb{S}_\mu\right|^2dxdt
  \leq \int_{\Omega}\Big[\frac{1}{2}\varrho_0|\mathbf{u}_0|^2
  +\frac{\varrho_0^\gamma}{\gamma - 1}\Big]dx
\end{equation}
for all $\tau\in[0,T]$, where $\nu$ is the same parameter appearing in
\eqref{constraint-nu}-\eqref{ass-3}.
\end{definition}



\begin{remark}
  A discussion concerning the derivation of \eqref{ee} is presented in
  the Appendix. In this regard, note that $\mu(\varrho)$ is a positive
  quantity, while $\lambda(\varrho)$ could be negative for
  $\varrho > 0$. Again, in obtaining \eqref{ee}, the viscous term
  involving $\lambda(\varrho)$ is directly controlled making use of
  the other viscous term, i.e.\ the one about $\mu(\varrho)$ by
  exploiting the elementary estimate
  $|\mathrm{div}\, \mathbf{u}|^2\leq 3
  |\mathbb{D}(\mathbf{u})|^2$. Through this last estimate, a minor
  part of information related to the considered system is lost,
  effectively determining an ``energy inequality". However, relation
  \eqref{ee} contains the most significant part of the energy
  contribution, and such an information is enough for our purposes.
\end{remark}

\noindent
The following existence result has been proved by Bresch, Desjardins
and G\'erard-Varet \cite{BDGV}.

\begin{theorem} \label{Th-ws} Let $T > 0$ be arbitrarily fixed.  Let
  also the initial data $(\varrho_0,\mathbf{u}_0)$ for the
  compressible Navier-Stokes system \eqref{cont}-\eqref{mom} be such
  that
  \begin{align*}
    \varrho_{0} \geq 0, \quad \,\varrho_{0} \in L^\gamma(\Omega), \quad\, 
    \sqrt{\varrho_0}\, \nabla s(\varrho_0)\in L^2(\Omega),\quad\,
    \varrho_{0}\mathbf{u}_{0} = 0 \,\,\,  \text{if} \,\,\, \varrho_{0}=0, \quad\,
    |\varrho_{0}\mathbf{u}_{0}|^2/\varrho_{0} \in L^1(\Omega).
  \end{align*}
  Then, for every $\varepsilon>0$ and $r_1>0$, there exists at least a
  weak solution $(\varrho, \mathbf{u})$ to
  system~\eqref{cont}-\eqref{mom} on the time interval $[0,T]$ with
  boundary conditions \eqref{bc}-\eqref{bc-bis}.
\end{theorem}

The existence of weak solutions for system \eqref{cont}-\eqref{mom}
without the drag term $r_1 \varrho \, |\mathbf{u}|\mathbf{u}$ has been shown
by Bresch,  Vasseur and Yu \cite{BVY} on the three-dimensional torus. In such
a paper, the authors extended the previous result by Vasseur and Yu \cite{VY}
to non-linear density dependent viscosity coefficients.

\subsection{Strong solution to the compressible Euler system}
\noindent
We recall the classical result concerning the local existence of
strong solutions for the compressible Euler system (see
\cite{Ag,Ve,Eb1,Eb2,Sc}).

\begin{theorem} \label{E-str} Let
  $(\varrho_0^E, \mathbf{u}_0^E) \in C^{1+\delta}(\Omega)$ with $\delta > 0$
  be some compatible initial data satisfying
  \begin{equation*}
    \inf_\Omega \varrho_0^E>0 \quad \,   \text{and}   \quad \,  \sup_\Omega \varrho_0^E < \infty.
  \end{equation*}
  Then, there exists a time $T^E>0$ depending on
$(\varrho_0^E, \mathbf{u}_0^E)$ and a unique strong solution of
\eqref{cont-E}-\eqref{bc-E} on the time interval $[0,T^E]$ satisfying
\begin{equation} \label{E-int} (\varrho^E, \mathbf{u}^E)
  \in C_w(0,T^E; C^{1+\delta}(\Omega)) \cap C^1([0,T^E]\times
  \overline{\Omega})
\end{equation}
and
\begin{equation} \label{E-vac} \inf_{[0,T^E]\times\Omega} \varrho^E >0
  \quad\, \text{and} \quad \, \sup_{[0,T^E]\times\Omega} \varrho^E <
  \infty.
\end{equation}
\end{theorem}

\noindent
In the statement above ``compatible'' refers to some
conditions satisfied by the initial data on the boundary
$\partial \Omega$ which are necessary for the existence of strong
solutions (see \cite{RaMa,Sc,Su}).

\section{Statement of the main results}
\noindent
We are now in a position to state the main results of the paper. From now on,
we suppose that the drag term parameter $r_1>0$ 
appearing in \eqref{mom} is a function of the viscosity parameter $\varepsilon>0$
such that $r_1=r_1(\varepsilon)\to 0$ as $\varepsilon\to 0$. 
This allows us to eliminate the dependence on $r_1$ of the weak solutions to 
the Navier-Stokes system \eqref{cont}-\eqref{mom}.

\begin{theorem}
\label{main-1}
Let
$(\varrho^E_0,\mathbf{u}^E_0)$ be as in Theorem~\ref{E-str} and
let $(\varrho^E,\mathbf{u}^E)$
be the corresponding strong solution of the compressible Euler system
\eqref{cont-E}-\eqref{bc-E} on the time interval $[0,T^E]$.
Let also $c>0$ and $0<T\leq T^E$ be arbitrarily fixed. 
Moreover, let 
$\varrho_0=\varrho_{0,\varepsilon}$ and $\mathbf{u}_{0}=\mathbf{u}_{0,\varepsilon}$ be
as in~Theorem \ref{Th-ws} and
assume that
\begin{equation} \label{id-conv}
\| \varrho_0 - \varrho^E_0 \|_{L^\gamma(\Omega)}
+ \int_\Omega \varrho_0 \left| \mathbf{u}_0 - \mathbf{u}^E_0 \right|^2
dx \to 0 \quad\, \text{as} \quad\, \varepsilon\to 0.
\end{equation}
Finally, let $\varrho =\varrho_{\varepsilon}$ and
$\mathbf{u}=\mathbf{u}_{\varepsilon}$ be a weak solution on the time
interval $[0,T]$ of the compressible Navier-Stokes
system~\eqref{cont}-\eqref{mom} with boundary conditions
\eqref{bc}-\eqref{bc-bis} corresponding to
$(\varrho_{0,\varepsilon},\mathbf{u}_{0,\varepsilon})$.  Assume that
\begin{equation} \label{k-conv}
 \int_{0}^{T} \int_{\Gamma_{c\varepsilon}} \left(\frac{\varrho^\gamma}{\varepsilon}
    +\varepsilon \, \frac{\varrho|\mathbf{u}|^2}{d_\Omega^2} \right) dxdt
  \to 0 \quad\,\textrm{as}\quad\, \varepsilon \to 0.
\end{equation}
Then, it follows that
\begin{equation} \label{tesi}
  \sup_{\tau \in [0,T]}  \left( \| \varrho - \varrho^E \|_{L^\gamma(\Omega)}
    + \|\varrho \mathbf{u} - \varrho^E\mathbf{u}^E\|_{L^1(\Omega)} \right)
  \to 0 \quad\, \text{as}  \quad\, \varepsilon\to 0.
\end{equation}
\end{theorem}

\begin{theorem} \label{main-2}
  Let $(\varrho^E_0,\mathbf{u}^E_0)$ be as in Theorem~\ref{E-str} and
  let $(\varrho^E,\mathbf{u}^E)$ be the corresponding strong solution
  of the compressible Euler system \eqref{cont-E}-\eqref{bc-E} on the
  time interval $[0,T^E]$.  Let also $c>0$ and $0<T\leq T^E$ be
  arbitrarily fixed.  Moreover, let
  $\varrho_0=\varrho_{0,\varepsilon}$ and
  $\mathbf{u}_{0}=\mathbf{u}_{0,\varepsilon}$ be as in~Theorem
  \ref{Th-ws} and assume that
  \begin{equation} \label{id-conv2}
    \| \varrho_0 - \varrho^E_0\|_{L^\gamma(\Omega)} + \int_\Omega \varrho_0 \left|
      \mathbf{u}_0 - \mathbf{u}^E_0 \right|^2 dx \to 0 \quad\,
    \text{as} \quad\, \varepsilon\to 0.
  \end{equation}
  Finally, let $\varrho =\varrho_{\varepsilon}$ and
  $\mathbf{u}=\mathbf{u}_{\varepsilon}$ be a weak solution on the time
  interval $[0,T]$ of the compressible Navier-Stokes
  system~\eqref{cont}-\eqref{mom} with boundary conditions
  \eqref{bc}-\eqref{bc-bis} corresponding to
  $(\varrho_{0,\varepsilon},\mathbf{u}_{0,\varepsilon})$.  Assume that
\begin{equation} \label{k-conv2}
\int_{0}^{T} \int_{\Gamma_{c\varepsilon}}
\left(\varepsilon \, \frac{\varrho|\mathbf{u}|^2}{d_\Omega^2}
+\varepsilon \, \frac{\varrho^2(\mathbf{u}\cdot \mathbf{n})^2}{d_\Omega^2}
+\varepsilon^{\left(\frac{\gamma-1}{\gamma}-\frac{1}{\gamma}\left(\frac{1-\nu}{3\nu}\right)\right)}
|\mathbb{S}_\mu|^2 \right) dxdt \to 0 \quad\, \textrm{as} \quad\, \varepsilon \to 0,
\end{equation}
where $\frac{\gamma-1}{\gamma}-\frac{1}{\gamma}\left(\frac{1-\nu}{3\nu}\right)= 1-\frac{1}{\gamma}\left(\frac{2}{3}
  + \frac{1}{3\nu}\right) $ {\rm (}see \eqref{constraint-nu} and \eqref{rho-split}{\rm )}, 
  and $\varrho^2(\mathbf{u}\cdot \mathbf{n})^2$ is
  meant as $\varrho(\sqrt{\varrho}\mathbf{u}\cdot \mathbf{n})^2$.
Then, it follows that
\begin{equation*}
  \sup_{\tau \in [0,T]} \left( \| \varrho - \varrho^E\|_{L^\gamma(\Omega)}
    + \|\varrho \mathbf{u} - \varrho^E\mathbf{u}^E \|_{L^1(\Omega)} \right)
  \to 0  \quad\, \text{as} \quad\, \varepsilon\to 0.
\end{equation*}
\end{theorem}

\begin{corollary} \label{cor}
Assume that $\mu(\varrho)=\varrho$ and $\lambda(\varrho)=0$. Let all the 
hypotheses of Theorem \ref{main-2} be valid except the third assumption in 
\eqref{k-conv2} which is replaced by
\begin{equation}
\label{k-conv-COR}
\int_{0}^{T} \int_{\Gamma_{c\varepsilon}}
\varepsilon^{\frac{\gamma-1}{\gamma}}
|\mathbb{S}_\mu|^2dxdt
\to 0 \quad\, \textrm{as} \quad\, \varepsilon \to 0.
\end{equation}
Then the conclusions of Theorem \ref{main-2} hold true.
\end{corollary}

The remaining of the paper is devoted to the proofs of Theorems \ref{main-1}-\ref{main-2}.
Corollary \ref{cor} does not follow directly from Theorem \ref{main-2}, but rather from its proof.
See Remark \ref{remcor} for details.

\begin{remark}
In the same spirit of \cite{BC}, one can replace the first two 
assumptions in \eqref{k-conv2} with the sole condition
\begin{equation} \label{bye-bye}
\int_{0}^{T} \int_{\Gamma_{c\varepsilon}}
\varepsilon^{\frac{\gamma-1}{\gamma}} 
\frac{\varrho|\mathbf{u}|^2}{d_\Omega^2}dxdt
\to 0 \,\, \textrm{ as } \,\, \varepsilon \to 0.
\end{equation}
See Remark \ref{bypass} for details.
\end{remark}


\section{Relative energy inequality}
\noindent
In this section, we derive a relative energy inequality satisfied by
the weak solution $(\varrho, \mathbf{u})$. Here and in the sequel, we
denote by
\begin{equation*}
  H(\varrho) = \frac{\varrho^\gamma}{\gamma - 1}.
\end{equation*}
For $\tau\in[0,T]$ we introduce the relative energy functional (see
e.g. \cite{FNJ-2012})
\begin{equation} \label{relative-energy}
  \mathcal{E}(\varrho, \mathbf{u}\,|\, r, \mathbf{U}) (\tau) =
  \int_{\Omega} \frac{1}{2} |\sqrt{\varrho}\mathbf{u}
  - \sqrt{\varrho}\,\mathbf{U}|^2  (\tau)  dx + \int_{\Omega}  H(\varrho|r)(\tau)  dx,
\end{equation}
where $(\varrho,\mathbf{u})$ is a weak solution to
\eqref{cont}-\eqref{mom} while $(r,\mathbf{U})$ is a pair of smooth
test functions such that $\mathbf{U}|_{\partial \Omega}=0$ and
$r>0$. Moreover
\begin{equation} \label{H-rho-r}
  H(\varrho|r) = H(\varrho) - H(r) - H'(r)(\varrho - r).
\end{equation}
Note that $H(\varrho|r)$ is non-negative and strictly convex. In
addition, it is equal to zero when $\varrho = r$ and it grows at
infinity as $\varrho^\gamma$.  Consequently, for any compact set
$K\subset (0,\infty)$ there exist two constants $c_3,c_4>0$ such that
for every $\varrho\geq0$ and every $r\in K$
\begin{equation}
  \label{puffo}
  c_3\big(\left|\varrho
    -r\right|^{2}1_{\left\{ \left|\varrho-r\right|<1\right\} }+\left|\varrho
    -r\right|^{\gamma}1_{\left\{ \left|\varrho-r\right|\geq 1\right\} }\big)
  \leq H(\varrho|r)
\end{equation}
and
\begin{equation}
  \label{puffo2}
  H(\varrho|r)\leq c_4\big(\left|\varrho
    -r\right|^{2}1_{\left\{ \left|\varrho-r\right|<1\right\} }+\left|\varrho
    -r\right|^{\gamma}1_{\left\{ \left|\varrho-r\right|\geq 1\right\} }\big)
\end{equation}
(see e.g.\ \cite{FNS-2011} and \cite[Section~2.1]{Su}).  Moreover,
being $\Omega$ a bounded domain, for any compact set
$K\subset (0,\infty)$ there exists a constant $c_5>0$ such that for
every $\varrho:\Omega\to [0,\infty)$ and $r:\Omega\to K$
\begin{equation} \label{H-gamma-1} \| \varrho - r
  \|_{L^\gamma(\Omega)}^\gamma \leq c_5\left(\int_\Omega
    H(\varrho|r)dx \right)^{\frac{\gamma}{2}} + c_5\int_\Omega
  H(\varrho|r)dx,
\end{equation}
and
\begin{equation} \label{H-gamma-2} \int_\Omega H(\varrho|r)dx \leq
  c_5\| \varrho - r \|_{L^\gamma(\Omega)}^{\gamma} +c_5\| \varrho - r
  \|_{L^\gamma(\Omega)}^2
\end{equation}
(see again \cite[Section~2.1]{Su}).
Thanks to the energy inequality \eqref{ee}, we infer that
\begin{align*}
\mathcal{E}(\varrho, \mathbf{u}\,|\, r, \mathbf{U}) (s)\big|_{s=0}^{s=\tau} 
&\leq  \int_\Omega \left(\frac{1}{2} \varrho |\mathbf{U}|^2 - \varrho \mathbf{u} 
\cdot \mathbf{U} \right)(s,\cdot) dx\big|_{s=0}^{s=\tau}
- \int_\Omega \left( H(r) + H'(r)(\varrho-r) \right) (s,\cdot)dx \big|_{s=0}^{s=\tau}\\
&\quad- \varepsilon\nu\int_0^\tau\int_{\Omega}\left|\mathbb{S}_\mu\right|^2dxdt.
\end{align*}
Next, we test the continuity equation by $\frac{1}{2}|\mathbf{U}|^2$ and the momentum equation by $\mathbf{U}$,
to get
\begin{equation*}
\int_0^\tau \int_{\Omega}\frac{1}{2}\varrho \partial_t |\mathbf{U}|^2\, dx dt 
+\int_0^\tau \int_{\Omega}\frac{1}{2} \varrho \mathbf{u} \cdot \nabla |\mathbf{U}|^2 \,dx dt  
= \int_\Omega \frac{1}{2} \varrho |\mathbf{U}|^2 (s,\cdot) dx\big|_{s=0}^{s=\tau}
\end{equation*}
and
\begin{gather*}
-\int_0^\tau \int_{\Omega} \varrho \mathbf{u}\cdot \partial_t \mathbf{U} \, dxdt
-\int_0^\tau \int_{\Omega} (\sqrt{\varrho} \mathbf{u} 
\otimes \sqrt{\varrho} \mathbf{u}) : \nabla \mathbf{U} \, dxdt 
-\int_0^\tau \int_{\Omega} p(\varrho) \textrm{div}\,\mathbf{U} \,dxdt \\\noalign{\vskip0.5mm}
+\varepsilon\int_0^\tau \int_{\Omega} \Big[2 \Big(\sqrt{\mu(\varrho)}\mathbb{S}_\mu 
+\frac{\lambda(\varrho)}{2\mu(\varrho)}\textrm{Tr}\big(\sqrt{\mu(\varrho)}\mathbb{S}_\mu\big)\mathbb{I} \Big)
: \nabla \mathbf{U} \Big] dxdt 
+r_1\int_{0}^\tau\int_{\Omega}|\sqrt{\varrho}\mathbf{u}|\sqrt{\varrho}\mathbf{u}\cdot \mathbf{U} dxdt\\\noalign{\vskip0.3mm}
= -\int_\Omega (\varrho \mathbf{u} \cdot \mathbf{U})(s,\cdot) dx\big|_{s=0}^{s=\tau}.
\end{gather*}
Then, testing the continuity equation by $H'(r)$ and noticing that
$H'(r)r - H(r) = p(r)$, after some straightforward computation, we have
\begin{equation*}  
  \int_\Omega \left( H(r) + H'(r)(\varrho-r) \right) (s,\cdot)dx \big|_{s=0}^{s=\tau}
= \int_0^\tau\!\! \int_\Omega  \left[  \partial_t  (H'(r))(\varrho - r) + \varrho
    \mathbf{u} \cdot \nabla (H'(r))  \right]  dxdt.
\end{equation*}
Therefore, we obtain
\begin{align*}
\mathcal{E}(\varrho, \mathbf{u}\,|\, r, \mathbf{U}) (s)\big|_{s=0}^{s=\tau} 
&\leq  \int_0^\tau \int_{\Omega}  \frac{1}{2}\varrho \partial_t |\mathbf{U}|^2\, dx dt + \int_0^\tau \int_{\Omega}
\frac{1}{2}\varrho \mathbf{u} \cdot \nabla |\mathbf{U}|^2 \,dx dt
\\&\quad - \int_0^\tau \int_{\Omega} \varrho \mathbf{u}\cdot \partial_t \mathbf{U} \, dxdt
-\int_0^\tau \int_{\Omega} (\sqrt{\varrho} \mathbf{u} \otimes\sqrt{\varrho} \mathbf{u}) : \nabla \mathbf{U} \, dxdt 
\\&\quad +\varepsilon\int_0^\tau \int_{\Omega} \Big[2 \Big(\sqrt{\mu(\varrho)}\mathbb{S}_\mu 
+\frac{\lambda(\varrho)}{2\mu(\varrho)}\textrm{Tr}\big(\sqrt{\mu(\varrho)}\mathbb{S}_\mu\big)\mathbb{I} \Big)
: \nabla \mathbf{U} \Big] dxdt 
\\&\quad 
- \int_0^\tau\!\! \int_\Omega  \left[  \partial_t  (H'(r))(\varrho - r) 
+ \varrho\mathbf{u} \cdot \nabla (H'(r)) +p(\varrho) \mbox{div} \mathbf{U} \right]  dxdt
\\&\quad +r_1\int_{0}^\tau\int_{\Omega}|\sqrt{\varrho}\mathbf{u}|\sqrt{\varrho}\mathbf{u}\cdot \mathbf{U}dxdt
-\varepsilon\nu\int_0^\tau\int_{\Omega}\left|\mathbb{S}_\mu\right|^2dxdt. 
\end{align*}
We now consider $(\varrho^E, \mathbf{u}^E)$ the strong solution of the Euler system
\eqref{cont-E}--\eqref{mom-E}.
Rearranging the relative entropy inequality above with $r = \varrho^E$, we obtain
\begin{gather} \label{step-4}
\begin{aligned}
\mathcal{E}(\varrho, \mathbf{u}\,|\, \varrho^E, \mathbf{U}) (s)\big|_{s=0}^{s=\tau} 
 &\leq \int_0^\tau \int_{\Omega}
\left[\partial_t \mathbf{U} \cdot (\varrho\mathbf{U} - \varrho\mathbf{u}) 
+ (\sqrt{\varrho}\mathbf{u} \cdot \nabla) \mathbf{U} \cdot (\sqrt{\varrho}\mathbf{U} - \sqrt{\varrho}\mathbf{u})\right]\,dx dt
\\ & \quad+\varepsilon\int_0^\tau \int_{\Omega} \Big[2 \Big(\sqrt{\mu(\varrho)}\mathbb{S}_\mu 
+\frac{\lambda(\varrho)}{2\mu(\varrho)}\textrm{Tr}\big(\sqrt{\mu(\varrho)}\mathbb{S}_\mu\big)\mathbb{I} \Big)
: \nabla \mathbf{U} \Big] dxdt 
 \\ &\quad - \int_0^\tau\!\! \int_\Omega  \left[  \partial_t  (H'(\varrho^E))(\varrho - \varrho^E) + \varrho
\mathbf{u} \cdot \nabla (H'(\varrho^E))  +p(\varrho) \mbox{div}\mathbf{U}\right]  dxdt
\\&\quad +r_1\int_{0}^\tau\int_{\Omega}|\sqrt{\varrho}\mathbf{u}|\sqrt{\varrho}\mathbf{u}\cdot \mathbf{U}dxdt
-\varepsilon\nu\int_0^\tau\int_{\Omega}\left|\mathbb{S}_\mu\right|^2dxdt,
\end{aligned}
\end{gather}
having set $(\sqrt{\varrho}\mathbf{u}\cdot \nabla)(\, \cdot\, ) =\sum_{i=1}^3(\sqrt{\varrho}\mathbf{u})_i\,\D_i(\, \cdot\, )$.
We now multiply \eqref{cont-E} by $H'(\varrho^E)$. Performing integration by parts, we find
\begin{equation*}
\int_0^\tau \int_\Omega\left[H'(\varrho^E)\partial_t \varrho^E 
- \varrho^E \nabla (H'(\varrho^E)) \cdot \mathbf{u}^E \right]dxdt= 0.
\end{equation*}
From $p(\varrho^E) = H'(\varrho^E)\varrho^E - H(\varrho^E)$,
we also infer that $\varrho^E \nabla (H'(\varrho^E)) = \nabla
p(\varrho^E)$.  Accordingly,
performing again integration by parts, we obtain
\begin{equation*}
  \int_0^\tau \int_\Omega H'(\varrho^E)\partial_t \varrho^E dxdt
  + \int_0^\tau \int_\Omega p(\varrho^E) \div \mathbf{u}^E dxdt=0.
\end{equation*}
Back to \eqref{step-4}, we have
\begin{align*} 
\mathcal{E}(\varrho, \mathbf{u}\,|\, \varrho^E, \mathbf{U}) (s)\big|_{s=0}^{s=\tau} &\leq \int_0^\tau \int_{\Omega}
\left[\partial_t \mathbf{U} \cdot (\varrho\mathbf{U} - \varrho\mathbf{u}) 
+ (\sqrt{\varrho}\mathbf{u} \cdot \nabla) \mathbf{U} \cdot (\sqrt{\varrho}\mathbf{U} - \sqrt{\varrho}\mathbf{u})\right]\,dx dt
\\ \nonumber & \quad+\varepsilon\int_0^\tau \int_{\Omega} \Big[2 \Big(\sqrt{\mu(\varrho)}\mathbb{S}_\mu 
+\frac{\lambda(\varrho)}{2\mu(\varrho)}\textrm{Tr}\big(\sqrt{\mu(\varrho)}\mathbb{S}_\mu\big)\mathbb{I} \Big)
: \nabla \mathbf{U} \Big] dxdt 
  \\ \nonumber &\quad - \int_0^\tau\!\! \int_\Omega  \left[  \partial_t  (H'(\varrho^E))(\varrho - \varrho^E) + \varrho
\mathbf{u} \cdot \nabla (H'(\varrho^E))  +p(\varrho) \mbox{div}\mathbf{U}\right]  dxdt
\\ \nonumber
&\quad +\int_0^\tau \int_\Omega  H'(\varrho^E)\partial_t \varrho^E dxdt + \int_0^\tau \int_\Omega p(\varrho^E) \div \mathbf{u}^E dxdt
\\ \nonumber &\quad +r_1\int_{0}^\tau\int_{\Omega}|\sqrt{\varrho}\mathbf{u}|\sqrt{\varrho}\mathbf{u}\cdot \mathbf{U}dxdt
-\varepsilon\nu\int_0^\tau\int_{\Omega}\left|\mathbb{S}_\mu\right|^2dxdt.
  \end{align*}
Now, from \eqref{mom-E} and the continuity equation (\ref{cont-E}), we easily deduce that
\begin{equation*}
\partial_t \mathbf{u}^E \cdot
  (\varrho\mathbf{U}- \varrho\mathbf{u}) +
(\sqrt{\varrho}\mathbf{u}^E \cdot \nabla) \mathbf{u}^E  \cdot
  (\sqrt{\varrho}\mathbf{U}-\sqrt{\varrho}\mathbf{u}) +
  \frac{1}{\varrho^E}\nabla p(\varrho^E) \cdot
  (\varrho\mathbf{U}-\varrho\mathbf{u}) =0.
\end{equation*}
Hence, we arrive at
\begin{equation} \label{step-6}
\begin{aligned}
\mathcal{E}(\varrho, \mathbf{u}\,|\, \varrho^E, \mathbf{U}) (s)\big|_{s=0}^{s=\tau} &\leq \int_0^\tau \int_{\Omega}
\left[\partial_t \mathbf{U} \cdot (\varrho\mathbf{U} - \varrho\mathbf{u}) 
+ (\sqrt{\varrho}\mathbf{u} \cdot \nabla) \mathbf{U} \cdot (\sqrt{\varrho}\mathbf{U} - \sqrt{\varrho}\mathbf{u})\right]\,dx dt\\
& \quad-\int_0^\tau \int_{\Omega}[\partial_t \mathbf{u}^E \cdot
  (\varrho\mathbf{U}- \varrho\mathbf{u}) +
(\sqrt{\varrho}\mathbf{u}^E \cdot \nabla) \mathbf{u}^E  \cdot
  (\sqrt{\varrho}\mathbf{U}-\sqrt{\varrho}\mathbf{u})] dx dt\\
 & \quad -\int_0^\tau \int_{\Omega}\frac{1}{\varrho^E}\nabla p(\varrho^E) \cdot
  (\varrho\mathbf{U}-\varrho\mathbf{u}) dx dt
\\ & \quad+\varepsilon\int_0^\tau \int_{\Omega} \Big[2 \Big(\sqrt{\mu(\varrho)}\mathbb{S}_\mu 
+\frac{\lambda(\varrho)}{2\mu(\varrho)}\textrm{Tr}\big(\sqrt{\mu(\varrho)}\mathbb{S}_\mu\big)\mathbb{I} \Big)
: \nabla \mathbf{U} \Big] dxdt 
  \\ &\quad - \int_0^\tau\!\! \int_\Omega  \left[  \partial_t  (H'(\varrho^E))(\varrho - \varrho^E) + \varrho
\mathbf{u} \cdot \nabla (H'(\varrho^E))  +p(\varrho) \mbox{div}\mathbf{U}\right]  dxdt
\\&\quad +\int_0^\tau \int_\Omega  H'(\varrho^E)\partial_t \varrho^E dxdt + \int_0^\tau \int_\Omega p(\varrho^E) \div \mathbf{u}^E dxdt
\\&\quad +r_1\int_{0}^\tau\int_{\Omega}|\sqrt{\varrho}\mathbf{u}|\sqrt{\varrho}\mathbf{u}\cdot \mathbf{U}dxdt
-\varepsilon\nu\int_0^\tau\int_{\Omega}\left|\mathbb{S}_\mu\right|^2dxdt.
  \end{aligned}
\end{equation}
Since using again (\ref{cont-E})
\begin{equation*}
  \partial_t (H'(\varrho^E)) 
= - p'(\varrho^E)\div \mathbf{u}^E - \nabla (H'(\varrho^E)) \cdot \mathbf{u}^E,
\end{equation*}
we see that
\begin{align*} 
&-\int_0^\tau \int_{\Omega}\frac{1}{\varrho^E}\nabla p(\varrho^E)
\cdot (\varrho\mathbf{U}-\varrho\mathbf{u}) dx dt-\int_0^\tau \int_\Omega 
 \partial_t  (H'(\varrho^E))(\varrho - \varrho^E)dx dt \\
&- \int_0^\tau \int_\Omega \varrho
\mathbf{u}\cdot \nabla (H'(\varrho^E)) dxdt +\int_0^\tau \int_\Omega H'(\varrho^E)\partial_t \varrho^E dxdt\\
& = \int_0^\tau \int_\Omega p'(\varrho^E)(\varrho -\varrho^E)\div \mathbf{u}^E dxdt+ \int_0^\tau \int_\Omega \nabla
(H'(\varrho^E))\cdot (\varrho\mathbf{u}^E - \varrho\mathbf{U}) dxdt.
\end{align*}
Inserting the identity above into \eqref{step-6}, we finally obtain
\begin{align}
\label{step-7}
\mathcal{E}(\varrho, \mathbf{u}\,|\, \varrho^E, \mathbf{U}) (s)\big|_{s=0}^{s=\tau}
&\leq \int_0^\tau \int_{\Omega}
\left[\partial_t \mathbf{U} \cdot (\varrho\mathbf{U} - \varrho\mathbf{u}) 
 + (\sqrt{\varrho}\mathbf{u} \cdot \nabla) \mathbf{U} \cdot (\sqrt{\varrho}\mathbf{U}
 - \sqrt{\varrho}\mathbf{u})\right]\,dx dt\\\nonumber
& \quad-\int_0^\tau \int_{\Omega}[\partial_t \mathbf{u}^E \cdot
  (\varrho\mathbf{U}- \varrho\mathbf{u}) +
(\sqrt{\varrho}\mathbf{u}^E \cdot \nabla) \mathbf{u}^E  \cdot
  (\sqrt{\varrho}\mathbf{U}-\sqrt{\varrho}\mathbf{u})] dx dt
\\\nonumber & \quad+\varepsilon\int_0^\tau \int_{\Omega} \Big[2 \Big(\sqrt{\mu(\varrho)}\mathbb{S}_\mu 
+\frac{\lambda(\varrho)}{2\mu(\varrho)}\textrm{Tr}\big(\sqrt{\mu(\varrho)}\mathbb{S}_\mu\big)\mathbb{I} \Big)
: \nabla \mathbf{U} \Big] dxdt \\\nonumber
&\quad +\int_0^\tau \int_\Omega \left[p(\varrho^E)\div \mathbf{u}^E
- p(\varrho) \div \mathbf{U}
+ p'(\varrho^E)(\varrho-\varrho^E)\div \mathbf{u}^E  \right] dxdt
\\\nonumber &\quad + \int_0^\tau \int_\Omega \nabla (H'(\varrho^E))\cdot (\varrho\mathbf{u}^E - \varrho\mathbf{U}) dxdt
\\\nonumber&\quad +r_1\int_{0}^\tau\int_{\Omega}|\sqrt{\varrho}\mathbf{u}|\sqrt{\varrho}\mathbf{u}\cdot \mathbf{U}dxdt
-\varepsilon\nu\int_0^\tau\int_{\Omega}\left|\mathbb{S}_\mu\right|^2dxdt.
\end{align}

\section{``Fake" boundary layer and uniform bounds}
\noindent
Following \cite{Ka,Su}, we introduce a Kato-type ``fake" boundary layer defined as
\begin{equation} \label{fl}
\mathbf{v}_{bl} = \mathbf{v}_{bl}(\varepsilon)=\xi \Big(\frac{d_{\Omega} (x)}{c\varepsilon} \Big) \mathbf{u}^E,
\end{equation}
where $c > 0$ is the constant in the statements of Theorems
\ref{main-1}-\ref{main-2}, while $\xi : [0,\infty)
\to [0,\infty)$ is a smooth cut-off function such that 
$\xi(0)=1$,  $\xi(r)=0$ for $r\geq1$ and satisfying
\begin{equation*}
\|\xi\|_{L^\infty(0,\infty)}<\infty, \qquad \|\xi'\|_{L^\infty(0,\infty)}<\infty.
\end{equation*}
Note that $\mathbf{v}_{bl} = \mathbf{u}^E$ on $[0,T] \times \partial \Omega$ and it has support in 
$[0,T] \times \Gamma_{c\varepsilon}$.
Choosing from now on
\begin{equation} \label{defU}
  \mathbf{U}=\mathbf{u}^E - \mathbf{v}_{bl}
\end{equation} 
and noticing that
\begin{equation*}
\int_0^\tau \int_\Omega  p(\varrho^E)\text{div}\,\mathbf{v}_{bl}\, dxdt=-
\int_0^\tau \int_\Omega \varrho^E \mathbf{v}_{bl} \cdot \nabla (H'(\varrho^E)) dxdt ,
\end{equation*}
we get from \eqref{step-7} the inequality
\begin{equation} \label{step-9}
  \mathcal{E}(\varrho, \mathbf{u}\,|\, \varrho^E, \mathbf{U}) (s)\big|_{s=0}^{s=\tau}
  +\varepsilon\nu\int_0^\tau\int_{\Omega}\left|\mathbb{S}_\mu\right|^2dxdt
\leq \mathcal{R}_1+ \ldots +\mathcal{R}_{7},
\end{equation}
where
\begin{align*} 
\mathcal{R}_1 &= \int_0^\tau \int_{\Omega}\Big(\big((\sqrt{\varrho}\mathbf{u}-\sqrt{\varrho}\mathbf{u}^E)\cdot\nabla\big)
\mathbf{u}^E\Big)\cdot(\sqrt{\varrho}\mathbf{U}-\sqrt{\varrho}\mathbf{u})\,dxdt,
\\
 \mathcal{R}_2 &=- \int_0^\tau\int_\Omega   \div \mathbf{u}^E\left(  p(\varrho)-p(\varrho^E)
 - p'(\varrho^E)(\varrho-\varrho^E)\right) dx dt,
\\
\mathcal{R}_3 &=- \int_0^\tau\int_\Omega (\varrho^E-\varrho) \mathbf{v}_{bl} \cdot \nabla (H'(\varrho^E)) dx dt,
\\
\mathcal{R}_4 &=\int_0^\tau\int_\Omega \left(\text{div} \mathbf{v}_{bl}\right)\left(p(\varrho)-p(\varrho^E)\right)dx dt,
\\
  \mathcal{R}_5 &= r_1 \int_0^\tau\int_{\Omega}|\sqrt{\varrho}\mathbf{u}|\sqrt{\varrho}\mathbf{u}\cdot \mathbf{U}dx dt,
  \\
\mathcal{R}_6 &= -\int_0^\tau\int_{\Omega}\left(\sqrt{\varrho}\partial_t \mathbf{v}_{bl}
+\left(\sqrt{\varrho}\mathbf{u}\cdot\nabla\right)\mathbf{v}_{bl}\right)\cdot(\sqrt{\varrho}\mathbf{U}-\sqrt{\varrho}\mathbf{u})dx dt,
\\
\mathcal{R}_7 &= \varepsilon \int_0^\tau\int_{\Omega} \Big[2 \Big(\sqrt{\mu(\varrho)}\mathbb{S}_\mu 
+\frac{\lambda(\varrho)}{2\mu(\varrho)}\textrm{Tr}\big(\sqrt{\mu(\varrho)}\mathbb{S}_\mu\big)\mathbb{I} \Big)
: \nabla \mathbf{U} \Big] dx dt.
\end{align*}
In the sequel, we denote by $C>0$ a generic structural constant
possibly depending on the initial data and $T>0$, but independent of
$\varepsilon$, and by $\eta(\varepsilon)$ a generic function
such that $\eta(\varepsilon)\to 0$ as $\varepsilon\to0$.
Our next goal is to show that under the assumptions of Theorems~\ref{main-1}
and \ref{main-2}, the terms $\mathcal{R}_i$, $i=1,\ldots,7$, can be controlled as follows
\begin{align} \label{goal}
&|\mathcal{R}_i| \leq C\int_0^\tau\mathcal{E}(\varrho, \mathbf{u}\,|\, \varrho^E, \mathbf{u}^E)(t) dt
+\eta(\varepsilon), \qquad i=1,\cdots 6,\\\label{goal2}
&|\mathcal{R}_7|
\leq \frac{\varepsilon\nu}{2} \int_0^\tau\int_{\Omega}\left|\mathbb{S}_\mu\right|^2dx dt
+ \eta(\varepsilon).
\end{align} 

\subsection{Uniform bounds on the weak solutions and the boundary layer}
\noindent
From \eqref{ee} and \eqref{id-conv}-\eqref{id-conv2}, it is readily seen that
for every $\varepsilon>0$ sufficiently small we have
\begin{equation} \label{uniform}
\| \varrho \|_{L^\infty(0,T;L^\gamma(\Omega))} \leq C, 
\qquad \| \sqrt{\varrho}\mathbf{u} \|_{L^\infty(0,T;L^2(\Omega))} \leq C,
\qquad \| \mathbb{S}_{\mu} \|_{L^2(0,T;L^2(\Omega))} \leq \frac{C}{\varepsilon^{1/2}}.
\end{equation} 
In addition, the ``fake" boundary layer $\mathbf{v}_{bl}$ constructed above satisfies
the following properties for every $\varepsilon>0$ sufficiently small (see \cite{BaNg,Su}):

\noindent
\begin{minipage}{0.4\linewidth}
\begin{align}
\label{prima}
 &\|\mathbf{v}_{bl}\|_{C([0,T]\times \overline\Omega)} \leq C,
  \\
\label{gradbl}   & \|\nabla \bv_{bl}\|_{L^\infty([0,T]\times\Omega)}\leq C \varepsilon^{-1}
  \\\label{dindon}
  &\|\text{div}\, \bv_{bl}\|_{C([0,T]\times \overline\Omega)}\leq C,
  \\
  & \|d_\Omega\nabla \bv_{bl}\|_{L^\infty(0,T; L^2(\Omega))}\leq C \varepsilon^{1/2},
\end{align}
\end{minipage}%
\begin{minipage}{0.6\linewidth}
\begin{align} 
\label{blLp} &\|\bv_{bl}\|_{C(0,T; L^p(\Omega))} \leq C
  \varepsilon^{1/p},\,\,\,  1\leq p <\infty, 
  \\\label{blLptime}
&\|\D_t \bv_{bl}\|_{C(0,T; L^p(\Omega))} \leq C
       \varepsilon^{1/p},\,\,\, 1\leq p <\infty,
  \\\label{dindondan}
  &\|\text{div}\, \bv_{bl}\|_{C(0,T; L^p(\Omega))} \leq C
    \varepsilon^{1/p},\,\,\, 1\leq p <\infty,
  \\\label{stimabella}
  &\|d^2_\Omega\nabla \bv_{bl}\|_{C([0,T]\times \overline\Omega)}
\leq C \varepsilon.
\end{align}
\end{minipage}

\medskip
\bigskip
\noindent
As a direct consequence of \eqref{defU} and \eqref{prima}, we have the further bound
\begin{equation}\label{boundU}
\|\mathbf{U}\|_{C([0,T]\times \overline\Omega)} \leq C.
\end{equation}

\section{Estimating the terms $\mathcal{R}_1$-$\mathcal{R}_5$} \label{secRi}
\noindent
The treatment of $\mathcal{R}_1$-$\mathcal{R}_4$ is basically the same
as in Sueur \cite{Su}, although here we have not used the same names
for these terms as in the latter reference.  Also, the term
$\mathcal{R}_5$ can be easily handled (see below). In what follows, it
is understood that we work with $\varepsilon>0$ sufficiently small in
a way that all the uniform bounds on the weak solutions and on the
boundary layer are valid.

\subsection{Estimate of $\mathcal{R}_1$}
\noindent
The term $\mathcal{R}_1$ can be decomposed as follows:
\begin{align*}
\mathcal{R}_1 = &- \int_0^\tau\int_{\Omega}\left(\left(
\sqrt{\varrho}\mathbf{u}-\sqrt{\varrho} \mathbf{u}^E\right)
\otimes \left( \sqrt{\varrho}\mathbf{u}-\sqrt{\varrho} \mathbf{u}^E\right)\right)
:\nabla \mathbf{u}^E\, dx dt\\
&- \int_0^\tau\int_{\Omega} \Big(((\sqrt{\varrho}\mathbf{u}-\sqrt{\varrho}\mathbf{u}^E
)\cdot\nabla)\mathbf{u}^E\Big)\cdot\sqrt{\varrho}\mathbf{v}_{bl}\, dx dt.
\end{align*}
Accordingly, thanks to \eqref{uniform} and \eqref{blLp}, we obtain
\begin{align*}
|\mathcal{R}_1| &\leq
C \int_0^\tau\int_{\Omega} |\sqrt{\varrho}\mathbf{u}-\sqrt{\varrho}\mathbf{u}^E|^2 dx dt
+ \left(\int_0^\tau\int_{\Omega}
|\sqrt{\varrho}\mathbf{u}-\sqrt{\varrho}\mathbf{u}^E|^2dx dt\right)^{1/2}\left(\int_0^\tau\int_{\Omega}
\varrho|\mathbf{v}_{bl}|^2dx dt\right)^{1/2}
\\
&\leq C \int_0^\tau\int_{\Omega} |\sqrt{\varrho}\mathbf{u}-\sqrt{\varrho}\mathbf{u}^E|^2 dx dt
+ C\left(\int_0^\tau\int_{\Omega}
|\sqrt{\varrho}\mathbf{u}-\sqrt{\varrho}\mathbf{u}^E|^2dx dt\right)^{1/2}\|\bv_{bl}\|_{C(0,T; L^{2\gamma/(\gamma-1)}(\Omega))}
\\
&\leq C \int_0^\tau\int_{\Omega} |\sqrt{\varrho}\mathbf{u}-\sqrt{\varrho}\mathbf{u}^E|^2 dx dt
     + C\left(\int_0^\tau\int_{\Omega}|\sqrt{\varrho}\mathbf{u}-\sqrt{\varrho}\mathbf{u}^E|^2\,dx dt\right)^{1/2}
     \varepsilon^{\frac{\gamma-1}{2\gamma}}.
\end{align*}
Exploiting Young's inequality, we finally get
\begin{equation} \label{est-1}
|\mathcal{R}_1| \leq C\int_0^\tau\int_{\Omega}
|\sqrt{\varrho}\mathbf{u}-\sqrt{\varrho}\mathbf{u}^E|^2dx dt +C\varepsilon^{\frac{\gamma-1}{\gamma}}
\leq C\int_0^\tau\mathcal{E}(\varrho, \mathbf{u}\,|\, \varrho^E, \mathbf{u}^E)(t) dt +\eta(\varepsilon).
\end{equation}

\subsection{Estimate of $\mathcal{R}_2$-$\mathcal{R}_4$}
\noindent
Now, consider the pressure terms $\mathcal{R}_2$-$\mathcal{R}_4$.
For $\mathcal{R}_2$, the following relation (see \eqref{H-rho-r})
holds true
\begin{align*}
|\mathcal{R}_2|&=\big|(\gamma-1)\int_0^\tau \int_\Omega 
( \div \mathbf{u}^E)H(\varrho|\varrho^E) dx dt \big|\\
&\leq C\int_0^\tau \int_\Omega  H(\varrho|\varrho^E) dx  dt\\ 
&\leq C\int_0^\tau\mathcal{E}(\varrho, \mathbf{u}\,|\, \varrho^E, \mathbf{u}^E)(t) dt. 
\end{align*}
Concerning $\mathcal{R}_3$, we decompose
\begin{align*}
\mathcal{R}_3 &=
-\int_0^\tau\int_{\Omega\cap\{|\varrho-\varrho^E|<1\}} (\varrho^E-\varrho) \mathbf{v}_{bl} \cdot \nabla (H'(\varrho^E))\, dx dt\\
&\quad -\int_0^\tau\int_{\Omega\cap\{|\varrho-\varrho^E|\geq1\}} (\varrho^E-\varrho) \mathbf{v}_{bl} \cdot \nabla (H'(\varrho^E))\, dx dt.
\end{align*}
Hence, invoking \eqref{puffo} and \eqref{blLp} together with the Young inequality,
we infer that
\begin{align*}
|\mathcal{R}_3| &\leq C
\int_0^\tau \int_\Omega  H(\varrho|\varrho^E) dx  dt 
+ C\|\bv_{bl}\|_{C(0,T; L^2(\Omega))}^2 + 
C\|\bv_{bl}\|_{C(0,T; L^{\gamma/(\gamma-1)}(\Omega))}^{\gamma/(\gamma-1)}\\
&\leq C\int_0^\tau \int_\Omega  H(\varrho|\varrho^E) dx  dt + \eta(\varepsilon)\\
&\leq C\int_0^\tau\mathcal{E}(\varrho, \mathbf{u}\,|\, \varrho^E, \mathbf{u}^E)(t) dt + \eta(\varepsilon). 
\end{align*}
Finally, we decompose $\mathcal{R}_4$ as
\begin{align*}
\mathcal{R}_4 &= \int_0^\tau \int_{\Omega}
\left(\text{div}\,\mathbf{v}_{bl}\right)
\left(p(\varrho)-p(\varrho^E)-p'(\varrho^E)(\varrho-\varrho^E)
\right)\, dx dt\\
&\quad + \int_0^\tau \int_{\Omega\cap\{|\varrho-\varrho^E|<1\}}
\left(\text{div}\,\mathbf{v}_{bl}\right)
p'(\varrho^E)(\varrho-\varrho^E) \, dx dt\\
&\quad + \int_0^\tau\int_{\Omega\cap\{|\varrho-\varrho^E|\geq1\}}
\left(\text{div}\,\mathbf{v}_{bl}\right)
p'(\varrho^E)(\varrho-\varrho^E) \, dx dt\\
&=(\gamma-1)\int_0^\tau \int_{\Omega}
\left(\text{div}\,\mathbf{v}_{bl}
\right)H(\varrho|\varrho^E) \, dx dt\\
&\quad + \int_0^\tau \int_{\Omega\cap\{|\varrho-\varrho^E|<1\}}
\left(\text{div}\,\mathbf{v}_{bl}
\right)p'(\varrho^E)(\varrho-\varrho^E) \, dx dt\\
&\quad + \int_0^\tau\int_{\Omega\cap\{|\varrho-\varrho^E|\geq1\}}
\left(\text{div}\,\mathbf{v}_{bl}
\right)p'(\varrho^E)(\varrho-\varrho^E) \, dx dt.
\end{align*}
Similarly as before, making use of \eqref{puffo}, \eqref{dindon},
\eqref{dindondan} and the Young inequality, we estimate
\begin{align*}
|\mathcal{R}_4| &\leq C\|\text{div}\, \bv_{bl}\|_{C([0,T]\times \overline\Omega)}
\int_0^\tau \int_\Omega  H(\varrho|\varrho^E) dx  dt +C\int_0^\tau \int_\Omega  H(\varrho|\varrho^E) dx  dt \\
&\quad + C\| \text{div}\,\bv_{bl}\|_{C(0,T; L^2(\Omega))}^2 + 
C\|\text{div}\,\bv_{bl}\|_{C(0,T; L^{\gamma/(\gamma-1)}(\Omega))}^{\gamma/(\gamma-1)}\\
&\leq C\int_0^\tau \int_\Omega  H(\varrho|\varrho^E) dx  dt + \eta(\varepsilon)\\
&\leq C\int_0^\tau\mathcal{E}(\varrho, \mathbf{u}\,|\, \varrho^E, \mathbf{u}^E)(t) dt + \eta(\varepsilon). 
\end{align*}

\subsection{Estimate of $\mathcal{R}_5$}
\noindent
Recalling that $r_1=r_1(\varepsilon)\to 0$ as $\varepsilon\to 0$,
and exploiting \eqref{uniform} together with \eqref{boundU},
the term~$\mathcal{R}_5$ can be handled easily as follows:
\begin{equation*}
  \mathcal{R}_5=r_1\int_{0}^\tau \int_{\Omega}|\sqrt{\varrho}\mathbf{u}|
  \sqrt{\varrho}\mathbf{u}\cdot \mathbf{U}\, dxdt 
\leq  \eta(\varepsilon).
\end{equation*}

\section{Proof of Theorem \ref{main-1}}
\noindent
Let us start with a preliminary observation. Under the
  assumption \eqref{k-conv}, we have in particular 
\begin{equation*}
\int_{0}^{T} \int_{\Gamma_{c\varepsilon}}
\frac{\varrho^\gamma}{\varepsilon} \, dxdt
\to 0 \quad\,\textrm{as}\quad\, \varepsilon \to 0,
\end{equation*}
and it follows immediately, for every $0< \ell \leq \gamma$, that
\begin{equation} \label{LUCA}
\int_{0}^{T} \int_{\Gamma_{c\varepsilon}}
\frac{\varrho^\ell}{\varepsilon}\, dxdt
\to 0 \quad\,\textrm{as}\quad\, \varepsilon \to 0.
\end{equation}
Indeed, since $|\Gamma_{c\varepsilon}|\leq C \varepsilon$, it is readily seen that
\begin{equation*}
  \int_{0}^{T} \int_{\Gamma_{c\varepsilon}} \frac{\varrho^\ell}{\varepsilon}\, dxdt
  \leq C \left(\int_{0}^{T} \int_{\Gamma_{c\varepsilon}}
\frac{\varrho^\gamma}{\varepsilon}\, dxdt\right)^{\frac{\ell}{\gamma}}
\to 0 \quad\,\textrm{as}\quad\, \varepsilon \to 0.
\end{equation*}

\subsection{Estimate of $\mathcal{R}_6$}
\noindent
As in \cite{Su}, we decompose the term $\mathcal{R}_6$ into the sum
\begin{align*}
\mathcal{R}_6 &= \int_0^\tau\int_{\Omega} \sqrt{\varrho}\partial_t \mathbf{v}_{bl}\cdot
\left(\sqrt{\varrho}\mathbf{u}-\sqrt{\varrho}\mathbf{u}^E \right)dx dt
+\int_0^\tau\int_{\Omega}\varrho \left(\partial_t\mathbf{v}_{bl}\cdot\mathbf{v}_{bl}\right)dx dt\\
&\quad -\int_0^\tau\int_{\Omega}
\big(\left(\sqrt{\varrho}\mathbf{u}\cdot\nabla\big)\mathbf{v}_{bl}\right)\cdot\sqrt{\varrho}\mathbf{U}\,dx dt
+\int_0^\tau\int_{\Omega} \left( \sqrt{\varrho}\mathbf{u}\otimes\sqrt{\varrho}\mathbf{u}\right)
: \nabla \mathbf{v}_{bl}\,dx dt\\\noalign{\vskip1mm}
&= \mathcal{R}_6^{(a)}+ \mathcal{R}_6^{(b)}+ \mathcal{R}_6^{(c)}+ \mathcal{R}_6^{(d)}.
\end{align*}
We have
\begin{equation*}
  |\mathcal{R}_6^{(a)}|
  \leq \left(\int_0^\tau\int_{\Omega}|\sqrt{\varrho}\mathbf{u}-\sqrt{\varrho}\mathbf{u}^E|^2dx dt
  \right)^{1/2} \left( \int_0^\tau\int_{\Omega}\varrho|\partial_t \mathbf{v}_{bl}|^2dx dt  \right)^{1/2}
\end{equation*}
and arguing as in the estimate of $\mathcal{R}_1$, but using
\eqref{blLptime} instead of \eqref{blLp}, we obtain
(see \eqref{relative-energy}) that
\begin{equation} \label{boundR6a}
|\mathcal{R}_6^{(a)}| \leq C\int_0^\tau\mathcal{E}(\varrho,
\mathbf{u}\,|\, \varrho^E, \mathbf{u}^E)(t) dt +\eta(\varepsilon).
\end{equation}
Moreover, in the light of \eqref{uniform}, together with \eqref{prima} and \eqref{blLptime},
it is readily seen that
\begin{equation} \label{boundR6b}
|\mathcal{R}_6^{(b)}|\leq C\|\partial_t \mathbf{v}_{bl}\|_{C(0,T; L^{\gamma/(\gamma-1)}(\Omega))}
\leq\eta(\varepsilon).
\end{equation}
Next, following \cite{Su}, we introduce the auxiliary scalar functions
\begin{equation*}
z(x)=\xi \Big(\frac{d_{\Omega} (x)}{c\varepsilon} \Big)\qquad \text{and}
\qquad \tilde z(x) = \frac{d_{\Omega}  (x)}{c\varepsilon}
  \xi' \Big(\frac{d_{\Omega} (x)}{c\varepsilon} \Big).
\end{equation*}
Accordingly, we can rewrite $\mathcal{R}_6^{(c)}$ in the form
\begin{equation*}
\mathcal{R}_6^{(c)}= - \int_0^\tau\int_{\Gamma_{c\varepsilon}}
    \sqrt{\varrho} z \mathbf{U}\cdot \left(\sqrt{\varrho}\mathbf{u}\cdot
\nabla\right)\mathbf{u}^E dx dt
+ \int_0^\tau\int_{\Gamma_{c\varepsilon}}
\frac{\varrho\mathbf{u}\cdot \mathbf{n}}{d_{\Omega}}
\tilde{z}\left(\mathbf{U}\cdot \mathbf{u}^E
\right)dx dt= \widetilde{\mathcal{R}}_6 + \widehat{\mathcal{R}}_6.
\end{equation*}
Being $\|z\|_{L^{\infty}(\Omega)}\leq C$, owing to \eqref{uniform}
and \eqref{boundU} we get
\begin{equation} \label{contr-tildeR6}
|\widetilde{\mathcal{R}}_6| \leq C \|z\|_{L^{2\gamma/(\gamma-1)}(\Gamma_{c\varepsilon})}
\leq C|\Gamma_{c\varepsilon}|^{\frac{\gamma-1}{2\gamma}} \leq \eta(\varepsilon).
\end{equation}
Moreover, since $\|\tilde z\|_{L^{\infty}(\Omega)}\leq C$, invoking
again \eqref{uniform} and \eqref{boundU} we find
\begin{equation*}
|\widehat{\mathcal{R}}_6| \leq C 
\left[\int_{0}^{T} \int_{\Gamma_{c\varepsilon}}
\frac{\varrho}{\varepsilon} dxdt\right]^{1/2}
\left[\int_{0}^{T} \int_{\Gamma_{c\varepsilon}}
\varepsilon \, \frac{\varrho|\mathbf{u}|^2}{d_\Omega^2}\, dx dt\right]^{1/2}.
\end{equation*}
The first term in the right-hand side above goes to zero as $\varepsilon\to0$
due to \eqref{LUCA}, while the second term goes to zero as $\varepsilon\to0$
due to our assumption \eqref{k-conv}. In summary
$|\widehat{\mathcal{R}}_6|\leq \eta(\varepsilon)$, yielding in turn
\begin{equation}
\label{boundR6c}
|\mathcal{R}_6^{(c)}| \leq\eta(\varepsilon).
\end{equation}
Finally, the term $\mathcal{R}_6^{(d)}$ can be handled as
\begin{align*}
|\mathcal{R}_6^{(d)}| &= \left|\int_0^\tau \int_{\Gamma_{c\varepsilon}}
\left(\frac{\sqrt{\varrho}\mathbf{u}}{d_\Omega}\otimes
\frac{\sqrt{\varrho}\mathbf{u}}{d_\Omega} \right):
d_\Omega^2\nabla \mathbf{v}_{bl}\, dx dt \right|\\\noalign{\vskip0.7mm}
&\leq \frac{C}{\varepsilon} \|d^2_\Omega\nabla \bv_{bl}\|_{C([0,T]\times \overline\Omega)}
\left[\int_{0}^{T} \int_{\Gamma_{c\varepsilon}}
\varepsilon \, \frac{\varrho|\mathbf{u}|^2}{d_\Omega^2}\, dx dt\right].
\end{align*}
In the light of \eqref{stimabella}, and
using once more our assumption \eqref{k-conv}, the right-hand side
above goes to zero as $\varepsilon\to0$, so that
\begin{equation}
\label{boundR6d}
|\mathcal{R}_6^{(d)}| \leq\eta(\varepsilon).
\end{equation}
Collecting \eqref{boundR6a}-\eqref{boundR6d} we conclude that
\begin{equation}
\label{boundR6}
|\mathcal{R}_6|\leq C\int_0^\tau\mathcal{E}(\varrho, \mathbf{u}\,|\,
\varrho^E, \mathbf{u}^E)(t) dt +\eta(\varepsilon).
\end{equation}

\subsection{Estimate of $\mathcal{R}_7$}
\noindent
Regarding the viscous term $\mathcal{R}_7$, thanks to \eqref{C-bound} we have
\begin{align*}
  |\mathcal{R}_7|&\leq C\varepsilon \int_0^\tau\int_{\Omega} \left|\sqrt{\mu(\varrho)}\mathbb{S}_\mu
  \right| \left| \nabla \mathbf{U} \right| dx dt
  \\
 &\leq C\varepsilon \int_0^\tau \int_{\Omega} \left|\sqrt{\mu(\varrho)}\mathbb{S}_\mu
 \right| \left| \nabla \mathbf{u}^E \right| dx dt
+ C\varepsilon \int_0^\tau \int_{\Omega} \left| \sqrt{\mu(\varrho)}\mathbb{S}_\mu
\right| \left| \nabla \mathbf{v}_{bl} \right| dx dt 
  = \widetilde{\mathcal{R}}_7 + \widehat{\mathcal{R}}_7.
\end{align*}
In light of \eqref{rho-split} and \eqref{constraint-nu} (which ensures that $2/3+1/{3\nu}\leq\gamma$),
together with \eqref{uniform}, we estimate
\begin{align*}
\widetilde{\mathcal{R}}_7
&\leq C \varepsilon \left[ \int_0^\tau\int_{\Omega} \mu(\varrho) dxdt \right]^{1/2}
\left[ \int_0^\tau\int_{\Omega}\left|\mathbb{S}_\mu\right|^2dx dt\right]^{1/2}
  \\
&=C \varepsilon\left[\int_0^\tau\int_{\Omega\cap\{\varrho< 1\}}\mu(\varrho)dx dt
+\int_0^\tau\int_{\Omega\cap\{\varrho\geq 1\}} \mu(\varrho)dx dt \right]^{1/2}\left[
\int_0^\tau \int_{\Omega}\left|\mathbb{S}_\mu\right|^2dx dt\right]^{1/2}
\\
&\leq C \varepsilon\left[\int_0^\tau\int_{\Omega\cap\{\varrho< 1\}}\varrho^{\left(\frac{2}{3}+\frac{\nu}{3}\right)}dx dt
+\int_0^\tau\int_{\Omega\cap\{\varrho\geq 1\}} \varrho^{\left(\frac{2}{3}+\frac{1}{3\nu}\right)}dx dt
\right]^{1/2}\left[ \int_0^\tau \int_{\Omega}\left|\mathbb{S}_\mu\right|^2dx dt\right]^{1/2}
\\
&\leq C \varepsilon\left[C+\int_0^\tau\int_{\Omega}\varrho^\gamma dx dt
\right]^{1/2} \left[ \int_0^\tau\int_{\Omega} \left| \mathbb{S}_\mu \right|^2dx dt\right]^{1/2}
\\
&\leq C \varepsilon\left[\int_0^\tau \int_{\Omega}\left| \mathbb{S}_\mu\right|^2dx dt \right]^{1/2}.
\end{align*}
Consequently
\begin{equation} \label{stimatilde}
\widetilde{\mathcal{R}} _7 \leq C \varepsilon \left[ \int_0^\tau \int_{\Omega}
\left|\mathbb{S}_\mu\right|^2dxdt\right]^{1/2}
\leq \frac{\varepsilon\nu}{2} \int_0^\tau\int_{\Omega}\left|\mathbb{S}_\mu\right|^2dx dt
+ \eta(\varepsilon),
\end{equation}
where the first term can be absorbed in the left-hand side of \eqref{step-9}.
Concerning the term $\widehat{\mathcal{R}}_7$, due to \eqref{gradbl} and
\eqref{uniform} we have
\begin{align*}
\widehat{\mathcal{R}}_7 &\leq C \varepsilon \|\nabla \mathbf{v}_{bl}
\|_{L^\infty([0,T]\times\Omega)} \left[\int_0^\tau \int_{\Gamma_{c\varepsilon}}\mu(\varrho)dxdt\right]^{1/2}
\left[\int_0^\tau \int_{\Gamma_{c\varepsilon}} \left|\mathbb{S}_\mu\right|^2 dxdt \right]^{1/2}
  \\
&\leq C\left[\int_0^\tau\int_{\Gamma_{c\varepsilon}}\mu(\varrho) dxdt\right]^{1/2}
\left[\int_0^\tau \int_{{\Omega}}\left|\mathbb{S}_\mu\right|^2 dxdt\right]^{1/2}
  \\
  &\leq C\left[\frac{1}{\varepsilon}\int_0^\tau \int_{\Gamma_{c\varepsilon}}
\mu(\varrho) dxdt \right]^{1/2}.
\end{align*}
As before, exploiting \eqref{rho-split}, we write
\begin{align*}
\frac{1}{\varepsilon}\int_0^\tau \int_{\Gamma_{c\varepsilon}} \mu(\varrho) dxdt
&= \frac{1}{\varepsilon}\int_0^\tau \int_{\Gamma_{c\varepsilon}\cap\{\varrho< 1\}}
\mu(\varrho) dxdt
+ \frac{1}{\varepsilon}\int_0^\tau \int_{\Gamma_{c\varepsilon}\cap\{\varrho\geq 1\}}
\mu(\varrho) dxdt\\\noalign{\vskip0.5mm}
&\leq
C \int_0^T \int_{\Gamma_{c\varepsilon}} \left[\frac{\varrho^{\left(\frac{2}{3}+\frac{\nu}{3}\right)}}{\varepsilon}\right]
dxdt +C \int_0^T \int_{\Gamma_{c\varepsilon}} \left[\frac{\varrho^{\left(\frac{2}{3}+\frac{1}{3\nu}\right)}}{\varepsilon}\right]dxdt.
\end{align*}
Since $0<2/3+\nu/3<\gamma$ and $0<2/3+1/{3\nu}\leq\gamma$, appealing to \eqref{LUCA} we infer that
$\widehat{\mathcal{R}}_7 \leq \eta(\varepsilon)$. Taking into account \eqref{stimatilde}, we end up with
\begin{equation} \label{stimaR7}
|\mathcal{R}_7| \leq \frac{\varepsilon\nu}{2} \int_0^\tau\int_{\Omega}\left|\mathbb{S}_\mu\right|^2dx dt
+ \eta(\varepsilon).
\end{equation}

\subsection{Conclusion of the proof of Theorem \ref{main-1}}
\noindent
So far, we have proved the validity of \eqref{goal}-\eqref{goal2}. Accordingly,
from \eqref{step-9}, we find
\begin{equation} \label{quasifinal}
\mathcal{E}(\varrho, \mathbf{u}\,|\, \varrho^E, \mathbf{U}) (\tau)\leq 
\mathcal{E}(\varrho, \mathbf{u}\,|\, \varrho^E, \mathbf{U}) (0)
+ C\int_0^\tau \mathcal{E}(\varrho, \mathbf{u}\,|\, \varrho^E, \mathbf{u}^E) (t)dt
+\eta(\varepsilon).
\end{equation}
Thanks to \eqref{uniform} and \eqref{blLp}, we have
\begin{align*}
\int_0^\tau \mathcal{E}(\varrho, \mathbf{u}\,|\, \varrho^E, \mathbf{u}^E) (t)dt
&\leq 
2\int_0^\tau \mathcal{E}(\varrho, \mathbf{u}\,|\, \varrho^E, \mathbf{U}) (t)dt
 + \int_0^\tau \int_\Omega \varrho |\mathbf{v}_{bl}|^2 dx dt
  \\
&\leq 2\int_0^\tau \mathcal{E}(\varrho, \mathbf{u}\,|\, \varrho^E, \mathbf{U}) (t)dt 
+C\|\bv_{bl}\|_{C(0,T; L^{2\gamma/(\gamma-1)}(\Omega))}^2
\\
&\leq 2\int_0^\tau \mathcal{E}(\varrho, \mathbf{u}\,|\, \varrho^E, \mathbf{U}) (t)dt 
+C \varepsilon^{\frac{\gamma-1}{\gamma}}.
\end{align*}
As a consequence, we infer from \eqref{quasifinal} that
\begin{equation*}
  \mathcal{E}(\varrho, \mathbf{u}\,|\, \varrho^E, \mathbf{U}) (\tau)
  \leq \mathcal{E}(\varrho, \mathbf{u}\,|\, \varrho^E, \mathbf{U}) (0)+\eta(\varepsilon)+ 
C\int_0^\tau \mathcal{E}(\varrho, \mathbf{u}\,|\, \varrho^E, \mathbf{U}) (t)dt,
\end{equation*}
and applying the Gronwall lemma we get
\begin{equation} \label{finale}
\sup_{\tau \in [0,T]}\mathcal{E}(\varrho, \mathbf{u}\,|\, \varrho^E,
\mathbf{U}) (\tau) \leq C \mathcal{E}(\varrho, \mathbf{u}\,|\, \varrho^E, \mathbf{U}) (0) + \eta(\epsilon).
\end{equation}
Next, invoking \eqref{H-gamma-2} and \eqref{blLp}, it is readily seen that
\begin{align*}
\mathcal{E}(\varrho, \mathbf{u}\,|\, \varrho^E, \mathbf{U}) (0)
&\leq C \int_\Omega \varrho_0 \left| \mathbf{u}_0 - \mathbf{u}^E_0
\right|^2 dx + C\int_\Omega \varrho_0 |\mathbf{v}_{bl}(0)|^2 dx
  \\
&\quad +C \|\varrho_0 - \varrho^E_0\|_{L^\gamma(\Omega)}^\gamma
+C \|\varrho_0 - \varrho^E_0\|_{L^\gamma(\Omega)}^2
  \\
&\leq C \int_\Omega\varrho_0\left|\mathbf{u}_0 - \mathbf{u}^E_0
\right|^2dx + C\|\varrho_0\|_{L^\gamma(\Omega)} \|\bv_{bl}\|_{C(0,T; L^{2\gamma/(\gamma-1)}(\Omega))}^2
  \\
&\quad +C \|\varrho_0 - \varrho^E_0\|_{L^\gamma(\Omega)}^\gamma
+C \|\varrho_0 - \varrho^E_0\|_{L^\gamma(\Omega)}^2
\\
&\leq C \int_\Omega\varrho_0\left|\mathbf{u}_0 - \mathbf{u}^E_0
\right|^2dx + C\varepsilon^{\frac{\gamma-1}{\gamma}}\|\varrho_0-\varrho_0^E\|_{L^\gamma(\Omega)} + 
 C\varepsilon^{\frac{\gamma-1}{\gamma}}\|\varrho_0^E\|_{L^\gamma(\Omega)}
  \\
&\quad +C \|\varrho_0 - \varrho^E_0\|_{L^\gamma(\Omega)}^\gamma+C \|
\varrho_0 - \varrho^E_0\|_{L^\gamma(\Omega)}^2.
\end{align*}
Therefore, due to \eqref{id-conv}, we have
$\mathcal{E}(\varrho, \mathbf{u}\,|\, \varrho^E, \mathbf{U}) (0)\to 0$
as $\varepsilon\to 0$. At this point, 
letting $\varepsilon\to 0$ in \eqref{finale}, we obtain
\begin{equation*}
\sup_{\tau \in [0,T]} \mathcal{E}(\varrho, \mathbf{u}\,|\, \varrho^E, \mathbf{U}) (\tau)\to0
\quad\, \text{as} \quad\, \varepsilon\to 0.
\end{equation*}
In the light of \eqref{H-gamma-1}, the latter yields
\begin{align} \label{rhofinal}
&\sup_{\tau \in [0,T]} \| \varrho - \varrho^E \|_{L^\gamma(\Omega)} \to 0
\quad\, \text{as} \quad\, \varepsilon\to 0,
  \\
\label{sqrtrhofinal}
&\sup_{\tau \in [0,T]}\|\sqrt{\varrho}
\mathbf{u} - \sqrt{\varrho}\mathbf{U}
\|_{L^2(\Omega)}\to 0
\quad\, \text{as} \quad\, \varepsilon\to 0.
\end{align}
Finally, noting that
\begin{equation*}
\varrho  \mathbf{u} - \varrho^E \mathbf{u}^E
=(\varrho  - \varrho^E)\mathbf{u}^E + \sqrt{\varrho}(\sqrt{\varrho}
\mathbf{u} - \sqrt{\varrho}\mathbf{U}) - \varrho \bv_{bl},
\end{equation*}
we estimate with the help of \eqref{uniform} and \eqref{blLp}
\begin{align*}
\|\varrho  \mathbf{u} - \varrho^E \mathbf{u}^E\|_{L^1(\Omega)}
&\leq C \| \varrho - \varrho^E \|_{L^\gamma(\Omega)} + 
C\|\sqrt{\varrho} \mathbf{u} - \sqrt{\varrho}\mathbf{U}\|_{L^2(\Omega)}
+ C\|\bv_{bl}\|_{C(0,T; L^{\gamma/(\gamma-1)}(\Omega))}
  \\
&\leq  C \| \varrho - \varrho^E \|_{L^\gamma(\Omega)} + 
C\|\sqrt{\varrho} \mathbf{u} - \sqrt{\varrho}\mathbf{U}
\|_{L^2(\Omega)} + C\varepsilon^{\frac{\gamma-1}{\gamma}}.
\end{align*}
Invoking \eqref{rhofinal}-\eqref{sqrtrhofinal},
we reach immediately the desired thesis \eqref{tesi}.
 \qed

\section{Proof of Theorem~\ref{main-2}}
\subsection{Estimate of $\mathcal{R}_6$}
\noindent
As in the proof of Theorem \ref{main-1}, we write
$\mathcal{R}_6=\mathcal{R}_6^{(a)}+ \mathcal{R}_6^{(b)}+
\mathcal{R}_6^{(c)}+ \mathcal{R}_6^{(d)}$.  The terms
$\mathcal{R}_6^{(a)}$ and $\mathcal{R}_6^{(b)}$ can be handled exactly
as before (see estimates \eqref{boundR6a}-\eqref{boundR6b}). In order
to control the term $\mathcal{R}_6^{(c)}$, we split it again into the
sum $\widetilde{\mathcal{R}}_6 + \widehat{\mathcal{R}}_6$. Noting that
$\widetilde{\mathcal{R}}_6$ satisfies \eqref{contr-tildeR6} we are left to deal with
the term $\widehat{\mathcal{R}}_6$. To this end, following \cite{Su} and
exploiting \eqref{boundU}, we estimate
\begin{equation*}
  |\widehat{\mathcal{R}}_6| \leq \frac{C}{\varepsilon^{1/2}}\|\tilde z\|_{L^{2}(\Gamma_{c\varepsilon})}
\left[\int_{0}^{T} \int_{\Gamma_{c\varepsilon}}
\varepsilon \, \frac{\varrho^2(\mathbf{u}\cdot \mathbf{n})^2}{d_\Omega^2}\, dx dt\right]^{1/2}.
\end{equation*}
It is readily seen that
\begin{equation*}
  \|\tilde z\|_{L^{2}(\Gamma_{c\varepsilon})}\leq 
  C|\Gamma_{c\varepsilon}|^{1/2}\leq C \varepsilon^{1/2}.
\end{equation*}
Thanks to our assumption \eqref{k-conv2}, we conclude that $|\widehat{\mathcal{R}}_6|\leq \eta(\varepsilon)$,
and we arrive at \eqref{boundR6c}. Finally, we observe that the remaining term $\mathcal{R}_6^{(d)}$ can be controlled
exactly as in the proof of Theorem \ref{main-1}, making use again of assumption \eqref{k-conv2}.
In summary, we have proved that \eqref{boundR6} holds true.

\subsection{Estimate of $\mathcal{R}_7$}
\noindent 
Arguing exactly as in the proof of Theorem \ref{main-1}, we have
$|\mathcal{R}_7|\leq \widetilde{\mathcal{R}}_7 + \widehat{\mathcal{R}}_7,$
where $\widetilde{\mathcal{R}}_7$ is controlled as in \eqref{stimatilde}.
Concerning $\widehat{\mathcal{R}}_7$,
an exploitation of \eqref{gradbl} yields
\begin{align*}
\widehat{\mathcal{R}}_7 &\leq C \varepsilon
\|\nabla \mathbf{v}_{bl} \|_{L^\infty([0,T]\times\Omega)}
\left[\int_0^\tau \int_{\Gamma_{c\varepsilon}}\mu(\varrho)dxdt\right]^{1/2}
\left[\int_0^\tau \int_{\Gamma_{c\varepsilon}}
\left|\mathbb{S}_\mu\right|^2 dxdt \right]^{1/2}
  \\
&\leq C\left[\int_0^\tau\int_{\Gamma_{c\varepsilon}}\mu(\varrho)dxdt\right]^{1/2}
\left[\int_0^\tau \int_{\Gamma_{c\varepsilon}}\left|\mathbb{S}_\mu\right|^2 dxdt\right]^{1/2}.
\end{align*}
Invoking \eqref{rho-split} and 
recalling that $|\Gamma_{c\varepsilon}|\leq C \varepsilon$, we obtain
\begin{align*}  
\int_0^\tau \int_{\Gamma_{c\varepsilon}} \mu(\varrho) dxdt
=\int_0^\tau \int_{\Gamma_{c\varepsilon}\cap\{\varrho< 1\}}
\mu(\varrho) dxdt
+\int_0^\tau \int_{\Gamma_{c\varepsilon}\cap\{\varrho\geq 1\}}
\mu(\varrho) dxdt\leq C\varepsilon + \int_0^\tau \int_{\Gamma_{c\varepsilon}}
\varrho^{\left(\frac{2}{3}+\frac{1}{3\nu}\right)} dxdt.
\end{align*}
Using again $|\Gamma_{c\varepsilon}|\leq C \varepsilon$, it follows
from \eqref{uniform} that
\begin{equation*}
\int_0^\tau \int_{\Gamma_{c\varepsilon}}
\varrho^{\left(\frac{2}{3}+\frac{1}{3\nu}\right)} dxdt
\leq C \varepsilon^{\left(1-\frac{1}{\gamma}\left(\frac{2}{3} + \frac{1}{3\nu}\right) \right)}.
\end{equation*}
Consequently (there is no harm in assuming $\varepsilon \leq 1$), we have
\begin{equation*}
\int_0^\tau \int_{\Gamma_{c\varepsilon}} \mu(\varrho)  dxdt
  \leq C \varepsilon^{\left(1-\frac{1}{\gamma}\left(\frac{2}{3} + \frac{1}{3\nu}\right) \right)}
  =C \varepsilon^{\left(\frac{\gamma-1}{\gamma}-\frac{1}{\gamma}\left(\frac{1-\nu}{3\nu}\right)\right)}.
\end{equation*}
Thanks to assumption \eqref{k-conv2},
we end up with 
\begin{equation*}
\widehat{\mathcal{R}}_7 \leq C \left[\int_0^T \int_{\Gamma_{c\varepsilon}}
\varepsilon^{\left(\frac{\gamma-1}{\gamma}-\frac{1}{\gamma}\left(\frac{1-\nu}{3\nu}\right)\right)}
\left|\mathbb{S}_\mu\right|^2 dxdt\right]^{1/2}
=\eta(\varepsilon).
\end{equation*}
In summary, estimate \eqref{stimaR7} holds true once again.

\subsection{Conclusion of the proof of Theorem~\ref{main-2}}
\noindent
As before, collecting all the bounds on the terms $\mathcal{R}_1$-$\mathcal{R}_7$ we obtain
\eqref{quasifinal}. The remaining part of the proof is identical to
the one of Theorem \ref{main-1}. \qed
\begin{remark} \label{remcor}
In order to prove Corollary \ref{cor}, one follows all the steps in the proof
of Theorem~\ref{main-2} except the estimate of the term 
$\widehat{\mathcal{R}}_7$ which, this time, can be controlled as
\begin{align*}
\widehat{\mathcal{R}}_7 &\leq C \varepsilon
\|\nabla \mathbf{v}_{bl} \|_{L^\infty([0,T]\times\Omega)}
\left[\int_0^\tau \int_{\Gamma_{c\varepsilon}}\varrho\, dxdt\right]^{1/2}
\left[\int_0^\tau \int_{\Gamma_{c\varepsilon}}
 \left|\mathbb{S}_\mu\right|^2 dxdt \right]^{1/2}
  \\
&\leq C \left[\int_0^T \int_{\Gamma_{c\varepsilon}}
\varepsilon^{\frac{\gamma-1}{\gamma}}\left|\mathbb{S}_\mu\right|^2 dxdt\right]^{1/2}.
\end{align*}
where the second inequality is a consequence of \eqref{gradbl}
together with $|\Gamma_{c\varepsilon}|\leq C \varepsilon$ and
\eqref{uniform}. Again, we infer that $\widehat{\mathcal{R}}_7\leq \eta(\varepsilon)$ due to assumption
\eqref{k-conv-COR}, and the thesis follows as before.
\end{remark}

\begin{remark}   \label{bypass}
  The first two assumptions in \eqref{k-conv2} have been used only to
  control $\widehat{\mathcal{R}}_6$ and $\mathcal{R}_6^{(d)}$. We now
  show that these two terms can be handled under the sole condition
  \eqref{bye-bye}. First we observe that, as a consequence of
  \eqref{uniform} and $|\Gamma_{c\varepsilon}|\leq C \varepsilon$, one  has
\begin{equation*}
\varepsilon^{\frac{1-\gamma}{\gamma}} 
\int_{0}^{T} \int_{\Gamma_{c\varepsilon}}
\varrho\, dxdt \leq C.
\end{equation*}
Accordingly, due to $\|\tilde z\|_{L^{\infty}(\Omega)}\leq C$
together with \eqref{uniform} and \eqref{boundU}, we infer that
\begin{equation*}
|\widehat{\mathcal{R}}_6| \leq C \left[\varepsilon^{\frac{1-\gamma}{\gamma}} 
\int_{0}^{T} \int_{\Gamma_{c\varepsilon}} \varrho\, dxdt\right]^{1/2}
\left[\int_{0}^{T} \int_{\Gamma_{c\varepsilon}} \varepsilon^{\frac{\gamma-1}{\gamma}} \,
  \frac{\varrho|\mathbf{u}|^2}{d_\Omega^2}\, dx dt\right]^{1/2}\leq \eta(\varepsilon),
\end{equation*}
where the second inequality follows from \eqref{bye-bye}. The term
$\mathcal{R}_6^{(d)}$ can be controlled exactly as in the proofs of
Theorems \ref{main-1}-\ref{main-2}:\ just note that \eqref{bye-bye} implies 
the first assumption in \eqref{k-conv2}.
\end{remark}

\subsection*{Acknowledgment}
\noindent
L. Bisconti and F. Dell'Oro are  members of the Gruppo Nazionale per
l'Analisi Mate\-ma\-tica, la Probabilit\`a e le loro Applicazioni (GNAMPA)
of the Istituto Nazionale di Alta Matematica (INdAM).
M. Caggio has been supported by the Praemium Academiae of \v S. Ne\v
casov\' a, and by the Czech Science Foundation under the grant GA\v CR 22-01591S.
The Institute of Mathematics CAS is supported by RVO:67985840.
This work has been written within the framework of the 
MUR grant Dipartimento di Eccellenza 2023-2027 of Dipartimento di Matematica - Politecnico di Milano.

\section*{Appendix}
\noindent 
We discuss briefly the derivation of the energy inequality \eqref{ee}.
To this end, let us consider smooth solutions $(\varrho,\mathbf{u})$
with $\varrho>0$. Multiplying the momentum equation \eqref{mom}
by $\mathbf{u}$ and after a standard computation, one can show that
\begin{equation} \label{eeAPP}
 \begin{aligned}
\int_{\Omega}\Big[\frac{1}{2}\varrho|\mathbf{u}|^{2}(\tau)
+H (\varrho)(\tau)\Big] dx
+ \varepsilon &\int_0^\tau\int_\Omega \big[2\mu(\varrho)|\mathbb{D}(\mathbf{u})|^2
+\lambda(\varrho)|\text{div}_x \mathbf{u}|^2|\big] dx dt
\\
&+ r_1 \int_0^\tau\int_\Omega \varrho |\mathbf{u}|^3 dx dt
 = \int_{\Omega}\Big[\frac{1}{2}\varrho_{0}|\mathbf{u}_{0}|^{2}+H(\varrho_{0})\Big]dx.
\end{aligned}
\end{equation}
We claim that
\begin{equation} \label{Lambda-intrappolato}
2\mu(\varrho)|\mathbb{D}(\mathbf{u})|^2
+\lambda(\varrho)|\text{div}_x \mathbf{u}|^2 \geq \nu \mu(\varrho)|\mathbb{D}(\mathbf{u})|^2
\end{equation}
where $\nu>0$ is the parameter appearing in assumption
\eqref{ass-1}-\eqref{ass-3}. This estimate is trivial when $\lambda(\varrho)\geq0$.
On the other hand, when $\lambda(\varrho)<0$, we use the
elementary inequality $|\text{div}_x \mathbf{u}|^2 \leq 3 |\mathbb{D}(\mathbf{u})|^2$.
Accordingly
\begin{equation*}
2\mu(\varrho)|\mathbb{D}(\mathbf{u})|^2
+\lambda(\varrho)|\text{div}_x \mathbf{u}|^2 
\geq [2\mu(\varrho) + 3\lambda(\varrho)]|\mathbb{D}(\mathbf{u})|^2
\end{equation*}
and the conclusion follow from \eqref{ass-3}.
Next, using \eqref{Lambda-intrappolato} into \eqref{eeAPP}, we obtain
\begin{equation} \label{eeNEW}
 \begin{aligned}
\int_{\Omega}\Big[\frac{1}{2}\varrho|\mathbf{u}|^{2}(\tau)
+H (\varrho)(\tau)\Big] dx
&+ \varepsilon\nu\int_0^\tau\int_\Omega \mu(\varrho)|\mathbb{D}(\mathbf{u})|^2 dx dt 
+ r_1\int_0^\tau\int_\Omega \varrho |\mathbf{u}|^3 dx dt
\\
&\leq \int_{\Omega}\Big[\frac{1}{2}\varrho_{0}|\mathbf{u}_{0}|^{2}+H(\varrho_{0})\Big]dx.
\end{aligned}
\end{equation}
Consider now a sequence of smooth solutions $(\varrho_n,\mathbf{u}_n)$
with $\varrho_n>0$ satisfying \eqref{eeNEW} and such that
\begin{equation*}
\int_{\Omega}\Big[\frac{1}{2}\varrho_{0n}|\mathbf{u}_{0n}|^{2}+H(\varrho_{0n})\Big]\leq E_0
\end{equation*}
with $E_0$ independent on $n$. 
From (\ref{eeNEW}), it follows that for all $\tau\in[0,T]$
\begin{align*}
&\|\sqrt{\varrho_n}\mathbf{u}_n\|_{L^\infty(0,\tau; L^2(\Omega))}\leq C(E_0),\\\noalign{\vskip0.7mm}
&\|\varrho_n\|_{L^\infty(0,\tau; L^\gamma(\Omega))}\leq C(E_0),\\\noalign{\vskip0.7mm}
&\|\sqrt{\mu(\varrho_n)}\mathbb{D}(\mathbf{u}_n)\|_{L^2((0,\tau)\times \Omega)}
 \leq C(\varepsilon,\nu,E_0),\\\noalign{\vskip0.7mm}
&\|\sqrt[3]{\varrho_n}\mathbf{u}_n\|_{L^3((0,\tau)\times \Omega)}\leq C(r_1,E_0).
\end{align*}
These bounds, combined with some further bounds coming from the
Bresch-Desjardins entropy inequality, have been used to prove the
stability of weak solutions and the existence result obtained in
\cite{BDGV} (see also \cite{MeVa} for the whole space and the periodic
case). In particular, the drag term provides the last bound in
${L^3((0,\tau)\times \Omega)}$ of the term
$\sqrt[3]{\varrho_n}\mathbf{u}_n$ which replaces the bound on
$\sqrt{\varrho_n}\mathbf{u}_n$ obtained in \cite{MeVa} by the use of
the multiplier $(1+\ln{(1+|\mathbf{u}|^2)})\mathbf{u}$ in the momentum
equation.  As is well-known, this control on the term
$\sqrt[3]{\varrho_n}\mathbf{u}_n$ is sufficient to obtain enough
compactness on $\sqrt{\varrho}_n\mathbf{u}_n$ in order to pass to the
limit.  Consequently, with the same arguments of \cite{BDGV,MeVa}, we
are able to pass to the limit as $n\to\infty$, getting the energy
inequality \eqref{ee}.

\end{document}